%

\magnification=\magstep1

\input amstex.tex
\long\def\comment#1\endcomment{}
\def\lable#1{{\rom [#1]}}
\def\nolable{\def\lable##1{}}

\nolable

\newcount\sectioncount
\newcount\commoncount\commoncount=1
\newcount\lemmacount\lemmacount=1
\newcount\defcount\defcount=1
\newcount\eqcount\eqcount=1
\newcount\theoremcount\theoremcount=1
\newcount\Ccount\Ccount=1

\def\firstsection#1{\sectioncount=#1 \advance\sectioncount by -1}
\firstsection1
\def\renewcount{\commoncount=1\eqcount=1
}
\renewcount

\def\newsection#1#2\par{\global\advance \sectioncount by 1%
\renewcount%
\specialhead\bigbf{\the\sectioncount}.\lable{#1}\bigbf\ #2
\expandafter\xdef\csname Section #1 \endcsname{%
\the\sectioncount}%
\endspecialhead
\relax}
\def\section#1{{\sl Section \csname Section #1 \endcsname}}

\def\newgrphno#1{
\the\sectioncount.\the\commoncount%
\expandafter\xdef\csname Grph #1 \endcsname{%
\the\sectioncount.\the\commoncount}%
\global\advance\commoncount by1\relax}
\def\newgrph#1{
\smallskip\noindent
{\bf \newgrphno{#1}.}\lable{#1}}%
\def\grph#1{{\csname Grph #1 \endcsname}}

\def\proof#1{
\smallskip\noindent{\bf\newgrphno{#1}.\lable{#1} Proof.}}

\def\newtheorem#1{
\proclaim{\the\sectioncount.\the\commoncount. Theorem}\lable{#1}
\expandafter\xdef\csname Theorem #1 \endcsname{\the\sectioncount.%
\the\commoncount}%
\global\advance\commoncount by1\relax}
\def\theorem#1{{\sl Theorem \csname Theorem #1 \endcsname}}

%

%

\def\newlemma#1{
\proclaim{\the\sectioncount.\the\commoncount.\lable{#1} Lemma}
\expandafter\xdef\csname Lemma #1 \endcsname{%
\the\sectioncount.\the\commoncount}%
\global\advance\commoncount by1\relax}
\def\lemma#1{{\sl Lemma \csname Lemma #1 \endcsname}}
\def\lemmano#1{\csname Lemma #1 \endcsname}

\def\newproposition#1{
\proclaim{\the\sectioncount.\the\commoncount.\lable{#1} Proposition}
\expandafter\xdef\csname Proposition #1 \endcsname{%
\the\sectioncount.\the\commoncount}%
\global\advance\commoncount by1\relax}
\def\proposition#1{{\sl Proposition \csname Proposition #1 \endcsname}}
\def\propositionno#1{\csname Proposition #1 \endcsname}

\def\newcorollary#1{
\proclaim{\the\sectioncount.\the\commoncount.\lable{#1} Corollary}
\expandafter\xdef\csname Corollary #1 \endcsname{%
\the\sectioncount.\the\commoncount}%
\global\advance\commoncount by1\relax}
\def\corollary#1{{\sl Corollary \csname Corollary #1 \endcsname}}
\def\corollaryno#1{\csname Corollary #1 \endcsname}

\def\newdefinition#1{
{\smallskip\noindent\bf\the\sectioncount.\the\commoncount.\lable{#1}
Definition.}
\expandafter\xdef\csname Definition #1 \endcsname{%
\the\sectioncount.\the\commoncount}%
\global\advance\commoncount by1\relax}
\def\definition#1{{\sl Definition \csname Definition #1 \endcsname}}
\def\definitionno#1{\csname Definition #1 \endcsname}
\let\dfntn=\definition

\def\makeeq#1{
\expandafter\xdef\csname Equation #1 \endcsname{%
\the\sectioncount.\the\eqcount}
\expandafter\xdef\csname Equationlable #1 \endcsname{{#1}}
\global\advance\eqcount by1\relax}
\def\eqqno#1{\csname Equation #1 \endcsname}
\def\eqq#1{{\rom (\eqqno{#1})}}
\def\neqq#1{\makeeq{#1} \eqno(\eqqno{#1})\rlap{\lable{#1}}}
\def\eqlable#1{{\rom{[\csname Equationlable #1 \endcsname]}}}

\documentstyle{amsppt}

\def\example#1{{\sl Example \csname Example #1 \endcsname}}
\def\qed{\ \ \hfill\hbox{$\square$}\smallskip}
\def\endstate{\endproclaim} 

\let\definition=\dfntn

\def\rom{\ifmmode \fam\rmfam \else\rm \fi}

\let\xpar=\par

\def\Month{\ifcase\month \or January\or February\or March\or April\or May\or
June\or July\or August\or September\or October\or November\or December\fi}

\font\bigbf=cmbx12

\thinmuskip = 2mu
\medmuskip = 2.5mu plus 1.5mu minus 2.1mu  
\thickmuskip = 4mu plus 6mu
\font\teneusm=eusm10
\font\seveneusm=eusm7
\font\fiveeusm=eusm5
\newfam\eusmfam
\textfont\eusmfam=\teneusm
\scriptfont\eusmfam=\seveneusm
\scriptscriptfont\eusmfam=\fiveeusm
\def\scr#1{{\fam\eusmfam\relax#1}}
\font\tenmib=cmmib10
\font\sevenmib=cmmib7
\font\fivemib=cmmib5
\newfam\mibfam
\textfont\mibfam=\tenmib
\scriptfont\mibfam=\sevenmib
\scriptscriptfont\mibfam=\fivemib
\def\mib{\fam\mibfam}
\font\tensf=cmss10
\font\sevensf=cmss8 scaled 833
\font\fivesf=cmr5
\newfam\sffam
\textfont\sffam=\tensf
\scriptfont\sffam=\sevensf
\scriptscriptfont\sffam=\fivesf
\def\sf{\fam\sffam}
\font\mathnine=cmmi9
\font\rmnine=cmr9
\font\cmsynine=cmsy9
\font\cmexnine=cmex10 scaled 913
\def\msmall#1{\hbox{$\displaystyle \font\ninesl=cmsl9
\textfont0=\rmnine \textfont1=\mathnine \textfont2=\cmsynine \
\textfont3=\cmexnine \textfont\slfam=\ninesl
{#1}$}}
\hyphenation{Lip-schit-zian Lip-schitz}
\def\cc{{\Bbb C}}

\def\nn{{\Bbb N}}

\def\zz{{\Bbb Z}}

\def\coker{{\sf Coker}\,}
\def\deg{{\sf deg}}
\def\dim{{\sf dim}\,}

\def\dimc{{\sf dim}_\cc}
\def\div{{\sf Div}}

\def\sfh{{\sf H}}
\def\hom{{\sf Hom}}

\def\id{{\sf Id}}

\def\inf{{\sf inf}\,}

\def\ker{{\sf Ker}\,}
\def\lim{\mathop{\sf lim}}
\def\log{{\sf log}}
\def\mod{{\sf mod}\,}

\def\supp{{\sf supp}\,}
\def\sup{{\sf sup}\,}
\def\sym{{\sf Sym}}

\def\epsi{\varepsilon}
\let\bs=\bss
\def\ogran{{\hskip0.7pt\vrule height7pt depth3pt\hskip0.7pt}}

\def\d{\partial}
\def\barr#1{\overline{#1}}
\def\dbar{{\overline\partial}}

\def\ddef{\mathrel{{=}\raise0.23pt\hbox{\rm:}}}
\def\deff{\mathrel{\raise0.23pt\hbox{\rm:}{=}}}
\def\ge{\geqslant}
\def\longto{\longrightarrow}
\def\inv{^{-1}}
\def\<{\langle}
\def\>{\rangle}
\def\fraction#1/#2{\mathchoice{{\msmall{ #1\over#2}}}%
{{ #1\over #2 }}{{#1/#2}}{{#1/#2}}}
\def\half{{\fraction1/2}}
\def\le{\leqslant}
\def\vph{^{\mathstrut}}
\def\lrar{\longrightarrow}

\def\maath{\mathsurround=0pt}
\def\lowminus{\hbox{\vbox to 0.3em{\vss\hbox{$ \maath\mkern-2mu \mathord- \mkern-2mu$}}}}

\def\deform #1#2#3#4#5{#1
\setbox1=\hbox{$\maath\scriptstyle\,\,#2\,$}
\buildrel{\,\,#2}\over{
\leftarrow\mkern-9mu
\hbox to\wd1{\cleaders\lowminus\hfill}}
#3
\setbox1=\hbox{$\maath\scriptstyle\,\,#4\,$}
\buildrel{#4\,\,}\over{
\hbox to\wd1{\cleaders\lowminus\hfill}
\mkern-9mu\rightarrow}
#5}

\def\mapdown#1|#2{\llap{$\vcenter{\hbox{$\scriptstyle #1$}}$}
 {{ \big\downarrow}}
  \rlap{$\vcenter{\hbox{$\scriptstyle #2$}}$}}

\def\emptyset{\varnothing}

\let\wh=\widehat
\let\wt=\widetilde
\let\ti=\tilde

\def\norm#1{\Vert{#1}\Vert}
\def\scirc{\mathchoice{\mathop{\msmall{\circ}}}{\mathop{\msmall{\circ}}}%
{{\scriptscriptstyle\circ}}{{\scriptscriptstyle\circ}}}
\def\state#1. {\medskip\noindent{\bf#1. }}
\def\qed{\ \hbox{ }\ \hbox{ }\ {\hfill$\square$}}
\def\Chi{\raise 2pt\hbox{$\chi$}}
\let\phI=\phi\let\phi=\varphi\let\varphi=\phI

\def\isig{{\mathchar"0106}} 
\def\iom{{\mathchar"010A}} 

\def\bfpsi{{\mathchar"0C20}}%
\def\bfphi{{\mathchar"0C27}}%

\def\eg{\hskip1pt plus1pt{\sl{e.g.\/\ \hskip1pt plus1pt}}}
\def\ie{\hskip1pt plus1pt{\sl i.e.\/\ \hskip1pt plus1pt}}
\def\iff{\hskip1pt plus1pt{\sl iff\/\hskip2pt plus1pt }}
\def\.{\thinspace}
\def\3{\ss}
\def\isl{\text{\sl i}}
\def\sli{{\sl i)} }              \def\slsi{{\sl S\hskip1pt i)} }
\def\slii{{\sl i$\!$i)} }        \def\slsii{{\sl S\hskip1pt i$\!$i)} }
\def\sliii{{\sl i$\!$i$\!$i)} }  \def\slsiii{{\sl S\hskip1pt i$\!$i$\!$i)} }
\def\sliv{{\sl i$\!$v)} }        \def\slsiv{{\sl S\hskip1pt i$\!$v)} }

\def\cal#1{{\scr{#1}}}
\def\cala{{\cal A}}
\def\alg(#1){\cala(\barr{#1})}
\def\calc{{\cal C}}
\def\cald{{\cal D}}

\def\calf{{\cal F}}
\def\cali{{\cal I}}

\def\calm{{\cal M}}
\def\caln{{\cal N}}
\def\calo{{\cal O}}

\def\calpi{{\cal{PI}}}

\def\calu{{\cal U}}

\def\Hom{{\cal H}{\mkern-3.5mu}{\italic o}{\mkern-3.2mu}{\italic m}}

\def\bff{{\mib F}}


\document
\baselineskip=12,8pt
\def\version{\hbox{\font\fiverm=cmr5 \fiverm Version of \the\day.\the\month.}}
\def\noversion{\def\version{}}
\noversion
\leftheadtext{\hss\vtop{%
\line{\hfil V.\.Shevchishin \hfil
\llap{\version}}%
\vskip 4pt \hrule }\hss}

\rightheadtext{\hss\vtop{%
\line{\rlap{\version}%
\hfil Moduli space of non-compact curves \hfil}%
\vskip 4pt \hrule }\hss}


\topmatter
\title \baselineskip=16pt \bigbf
A moduli space of non-compact curves
on a complex surface
\endtitle
\author
Vsevolod V.~Shevchishin
\endauthor
\abstract
We show that under mild boundary conditions the moduli space of
non-compact curves on a complex surface is (locally) an analytic subset of
a ball in Banach manifold, defined by {\sl finitely} 
many holomorphic function.
\endabstract
\address
Fakult\"at f\"ur Mathematik,
Ruhr-Universit\"at\.Bochum,\xpar
Universit\"atsstra\3e\.150,
D-44780~Bochum,~Germany.\xpar
\it Permanent Address: \rm
Institute of Applied Problem of Mechanics and Mathematics, 
Ukrainian Academy of Science, Lviv, Ukraine.
\endaddress
\email sewa\@cplx.ruhr-uni-bochum.de  \endemail
\thanks The author was supported by {\sl Graduiertenkolleg}
"Geometrie und Mathematische Physik" of Ruhr-University, Bochum.
\endthanks

\endtopmatter

\medskip
\smallskip
\centerline{September 1996.}

\smallskip
\centerline{Reviewed: Juli 1998.}

\bigskip
\bigskip
\firstsection0
\newsection{S0}{Introduction and basic notation.}

\nobreak
\newgrph{0.1}
In this paper $X$ denotes a {\sl smooth complex surface\/}. {\sl A curve $C$
in $X$} is an effec\-tive divisor. This means that $C$ is a locally finite
formal sum $C=\sum_i m_i C_i$, where every $C_i$ is a (closed) irreducible
analytic set of (co)dimension $1$, and $m_i$ are positive integers. We call
$C_i$ {\sl (irreducible) components of $C$ {\rm and} $m_i$ multiplicities}.
The set $|C| \deff \cup_i C_i$ is called the {\sl support of $C$}.
For an open  subset $U\subset X$ we define the {\sl restriction of $C$ to
$U$\/} as $C\cap U\deff \sum_i m_i (C_i\cap U)$.

With any component $C_i$ we associate the {\sl ideal sheaf\/}
$\cali_{C_i}$ whose group of sections over an open set $V\subset X$ is
$\Gamma(V,\cali_{C_i}) \deff \{\, f\in \Gamma(V,\calo)\;:\;
f\ogran_{V\cap C_i}\equiv 0\,\}$.
This is a coherent analytic sheaf and $\supp (\calo/\cali_{C_i})= C_i$.
Call $\cali_C \deff \prod_{i=1}^N \cali_{C_i}^{m_i}$
the {\sl ideal sheaf of $C$}, and $\calo_C \deff \calo_X/\cali_C$
the {\sl structure sheaf of $C$}.  The ideal sheaf $\cali_C$ is locally 
principle, \ie has locally the form $\cali_C\ogran_U= f_U \cdot \calo_X$.
We call such $f_U$ a {\sl local determining function} of $C$ in $U$ and $C\cap
U$ a {\sl divisor of $f_U$}, $C\cap U = \div(f_U)$.

The pair $(|C|, \calo_C)$ is a complex
subspace of $X$ (in general not reduced and reducible) which we shall denote
also by $C$. This means that we can concider $C$ as an {\sl analytic cycle} 
$C=\sum_i m_i C_i$ as also as {\sl subspace $C= (|C|, \calo_C)$} of $X$.

\newgrph{curr}
It is known (see \eg [Ha] or [Ki]), that one can associate to every curve 
$C=\sum_i m_i C_i$ a
{\sl closed positive integer (1,1)-current $\eta_C$} such that, for any 
continuous 2-form $\phi$ with the compact support in $X$,
$$
\eta_C(\phi)=\<\eta_C,\phi\>\deff \sum_i m_i\int_{C_i} \phi.
$$
Moreover, $C$ is completely determined by $\eta_C$ and every closed positive
integer (1,1)-current $\eta$ in $X$ corresponds to some curve, see [Ha] or
[Ki]. Thus we identify the set of curves in $X$ with a space 
$\calpi^{1,1}(X)$ of closed positive integer
(1,1)-currents in $X$ and induce the topology on a set of curves. Note
that $\calpi^{1,1}(X)$ is closed subset in the space $\cald'_2(X)$ of
2-currents in $X$.

The (weak) topology in the space $\calpi^{1,1}(X)$ gives us the
notion of weakly continuous family of curves in $X$. Namely, a family
$\{C_y\}_{y\in Y}$, parameterised by a topological space $Y$, is called {\sl
weakly continuous}, \iff the induced map $F:Y\to \calpi^{1,1}(X)$, $F(y)\deff
\eta\vph_{C_y}$, is continuous. We shall show later that this is equivalent
to the following condition: There exist an open covering $\{V_\alpha\}$ of
$X\times Y$ and continuous functions $f_\alpha \in C(V_\alpha, \cc)$ such
that for any $y\in Y$ the restriction of $f_\alpha$ on $(X\times\{y\})\cap
V_\alpha$ is holomorphic and generates the ideal sheaf $\cali_{C_y}=f_\alpha
\cdot \calo_{X\times \{y\}}$ of $C_y$. The $f_\alpha$ are called {\sl local
determining functions of a family $\{C_y\}$}.

In particular, a sequence $\{C_\nu\}$ of curves in $X$ {\sl converges weakly
to a curve $C_\infty$} \iff for any $x\in X$ there exist a neighbourhood
$V\ni x$ and a sequence of holomorphic functions $f_\nu\in \Gamma(V,\calo)$
which are determining for $C_\nu$ and which converge uniformly in $V$ to a
determining function $f_\infty$ of $C_\infty$.

\newgrph{0.2}
For the definition of the category of {\sl Banach analytic spaces}  we refer 
to [Dou], {\sl Section 3}. We note that for every Banach analytic space $Y$ 
and Banach spaces $E$ the sheaf $\calo_Y(E): U \subset Y \mapsto \Gamma(U, 
\calo_Y(E))$ of holomorphic $E$-valued morphisms between open subsets $U 
\subset Y$ and $E$ is a part of a definition of the structure of 
a Banach analytic space. In the case $E=\cc$ we denote this sheaf by 
$\calo_Y$. Any morphism $\bff: Y \to Z$ between two  Banach 
analytic spaces defines a continuous map $F: Y \to Z$ between corresponding
topological spaces, and a morphism of sheaves $F_E^\sharp : \calo_Y(E) \to F^*
\calo_Z(E)$ for any Banach space $E$. Here $F^* \calo_Z(E)$ denotes the 
pull-back of the sheaf $\calo_Z(E)$ w.r.t.\ continuous map $F$. Moreover,
a morphism $\bff :Y \to Z$ is defined by the data $F$ and $F_{(\cdot)}^\sharp$.

We say that a continuous map $F: Y \to Z$ is {\sl holomorphic} if it is induced
by a morphism $\bff: Y \to Z$. Note that such a morphism $\bff: Y \to Z$ 
can be not unique at the sheaf level. In particular, two different morphisms 
$\bff_1, \bff_2 \in {\rom Mor}(Y,E) = \Gamma(Y, \calo_Y(E))$ can induce the
same continuous map $F_1=F_2 : Y \to E$. This reflects the fact that a generic
Banach analytic space $Y$ is highly non-reduced.

\newdefinition{fin-codim}
We say that a Banach analytic space $Y$ is {\sl of finite type} \iff $Y$ can 
covered by local charts $Y_\alpha$ such that any $Y_\alpha$ is isomorphic to 
a zero set of a holomorphic map $f_\alpha:B_\alpha\to \cc^{n_\alpha}$, where
$B_\alpha$ denotes a ball in some Banach space. In particular, we have an
isomorphism
$
\calo_Y\ogran_{Y_\alpha} \cong
\calo_{B_\alpha}/(f_{\alpha,1},\cdots f_{\alpha,n_\alpha}),
$
where $(f_{\alpha,1},\cdots f_{\alpha,n_\alpha})$ denotes the ideal sheaf
generated by the components of $f_\alpha$. Such spaces are also refered to as
{\sl Banach analytic spaces of finite definition} or
{\sl Banach analytic spaces of finite codimension}.

\newdefinition{0.2a}
A {\sl holomorphic family} $\calc=\{C_y\}_{y\in Y}$ of curves in $X$ 
parametrised by a Banach analytic space $Y$ is given by an open covering
$\{V_\alpha\}$ of $X\times Y$ and holomorphic functions $f_\alpha \in
\Gamma(V_\alpha, \calo_{X\times Y})$ such that:

\sli if $V_\alpha \cap V_\beta \not = \emptyset$, then $f_\alpha= f_{\alpha
\beta} \cdot f_\beta$ for some {\sl invertible} $f_{\alpha\beta} \in \Gamma( 
V_\alpha \cap V_\beta, \calo^*_{X\times Y})$;

\slii for any $y\in Y$ the restriction of $f_\alpha$ on $V_\alpha \cap 
X\times \{y\}$ is not identically zero an is a local determining function
for a curve $C_y$.

The functions $f_\alpha$ are called {\sl local determining functions of
the family $\{C_y\}_{y\in Y}$}. 
The collection $\{ f_\alpha \}$ defines the sheaf of ideals $\cali_\calc
\subset \calo_{X\times Y}$ with $\cali_\calc\ogran_{V_\alpha} \deff f_\alpha
\cdot \calo_{X\times Y}$. Two holomorphic families parametrised by the same 
Banach analytic space $Y$ are {\sl isomorphic} \iff they define the same 
sheaf of ideals over $X\times Y$.

\smallskip
\newgrph{MainTh}
Now let $C^*$ be any curve in $X$ and $K\Subset |C^*|$ any compact subset
of its support. Our main result is the following

\smallskip
\proclaim{Main Theorem}
There exists an open set $U\subset X$ containing $K$ such that the set of
curves in $U$, which satisfy an appropriate boundary condition and are
sufficiently close to $C^*\cap U$, is a holomorphic family $\calc= \{\calc_t\}
_{t\in\calm}$ parameterised by a Banach analytic space $\calm$ of finite
type. Moreover, for every continuous (resp.\ holomorphic) family $\{C_y\}
_{y\in Y}$ with $C^*=C_{y_0}$ for some $y_0\in Y\!\!,$ there exist a 
neighbourhood $Y_0\subset Y$ of $y_0$ and a continuous map (resp.\ a morphism)
$F:Y_0\to \calm$ such that $C_y\cap U= \calc_{F(y)}$. Two such families
$\{C'_y\}\vph_{y\in Y}$ and $\{C''_y\}\vph_{y\in Y}$ coincide over $Y_0$ \iff
they induce the same continuous map (resp.\ morphism) $F: Y_0 \to \calm$.
\endstate

\newgrph{Corlls}
The theorem has several corollaries which are mainly due to the fact that
Banach analytic sets of finite type have sufficiently simple structure. In
particular, if $X$, $C^*$, and $U$ are as in the {\sl Main Theorem}, and if
$\{C_n\}$ is a sequence of curves in $X$ converging to $C^*$, then for any
$n>\!>1$ there exists a holomorphic family $\{C_\lambda\} _{\lambda \in
\Delta}$ of curves in $U$, which is parametrised by a disk $\Delta \subset
\cc$ and contains both $C_n\cap U$ and $C^*\cap U$. This allows to obtain
a generalization of the continuity principle of E.\.E.\.Levi.

\medskip
\newgrph{0.4}
The conclusion of the {\sl Main Theorem} is obtained by an explicit 
construction of the space $\calm$. The problem of deformation of a curve $C$ 
leads to study of the {\sl normal sheaf\/ $\caln_C$ to $C$} in $X$. It is 
defined as $\caln_C \deff \Hom_{\calo_X}\vph(\cali_C/\cali_C^2, \calo_C)$. 

To obtain a parametrising space as an analytic set in {\sl Banach} manifold
we introduce the notion of a (Banach) {\sl smoothness} $S$. This generalises
the usual smoothness classes such as $k$ times continuous differentiability
$C^k$, Sobolev smoothness $L^{k,p}$, or H\"older smoothness $C^{k,\alpha}$.
For such a smoothness $S$ we define a Banach space $\Gamma_S(C, \caln_C)\equiv
\sfh^0_S(C,\caln_C)$ of (holomorphic) sections of $\caln_C$ which are {\sl
$S$-smooth up to boundary $\d C$} (or simply $S$-smooth). 

The description of a moduli space $\calm$ in a neighbourhood of a marked 
point $y_0$ is usual for a deformation theory:

\newgrph{idea}
\sl There exists a ball $B\subset \sfh^0_S(C,\caln_C)$ and a holomorphic map
$\Phi:B \to \sfh^1(C,\caln_C)$ such that
\item{}\sli $\Phi(0)=0$, $d\Phi(0)=0$;
\item{}\slii $\Phi: B \to \sfh^1(C,\caln_C)$ is a local chart for $\calm$.
\item{}\sliii $0\in B$ corresponds to $y_0\in \calm$, parameterizing $C^*\cap
U$.
\rm

\noindent
The desired property of $\calm$ is based on the fact that non-compact
components of $C$ are Stein spaces, and consequently $\sfh^1(C, \caln_C)$
is finite-dimensional. In particular, 
$\calo_\calm\cong \calo_B/(\Phi_1,\ldots,\Phi_k),\qquad
k\deff \dim \sfh^1(C,\caln_C)$. Here $\Phi_i$ denote the components of $\Phi$
and $(\Phi_1,\ldots,\Phi_k)$ the ideal sheaf genereated by $\Phi_1,\ldots,
\Phi_k$.

\medskip
\newgrph{thanks} {\bf Acknoledgements.} I expresses his gratitude to
S.\.Ivashkovich for the constant support during preparation of this paper
and the pointing out possible applications. I also wish to thank 
D.\.Barlet, H.\.Flenner,
S.\.Kosarew, G.\.Schumacher, and B.\.Siebert for numerous helpful discussions
on different questions of the deformation theory.

\bigskip\bigskip
\newsection{Loc}{The local situation.}

\nobreak
\newgrph{1.1}
The construction of the moduli space $\calm$ is based on two special cases. 
One of them describes local deformations of curves and the other one allows 
to match two different local descriptions.

We first consider a local situation. For this we suppose that $\isig$ is
a smooth
complex curve with a smooth nonempty boundary $\d\isig$. Set $V \deff 
\isig \times \Delta$. 
Denote by $\Gamma_{L^\infty}(\isig,\calo^n)$ the Banach space of
$n$-tuples of holomorphic uniformly bounded functions on $\isig$. For every
such $f= (f_1(z),\ldots,f_n(z))\in \Gamma_{L^\infty} (\isig, \calo^n)$ we
define a {\sl Weierstra\3 polynomial}
$$
P_f(z,w)\deff w^n + \sum_{i=1}^n f_i(z)\,w^{n-i},
\qquad z\in \isig, \quad w\in \Delta
\neqq{WP}
$$
and a curve $C_f\subset V$ to be the zero divisor of $P_f(z,w)$.

\newlemma{1.2}
Every curve $C\subset V=\isig \times \Delta$ satisfying condition
$$
|C|\subset \isig\times \Delta(r)  \qquad \text{for some $r<1$}
\neqq{1.2}
$$
is a zero divisor of a uniquely defined Weierstra\3 polynomial $P_f(z,w)=w^n
+ \sum_{i=1}^n f_i(z)\,w^{n-i}$ with $f= (f_1,\ldots,f_n)\in
\Gamma_{L^\infty} (\isig, \calo^n)$.

The set $\calm^{(n)}_{L^\infty}(V)$ of those $f\in \Gamma_{L^\infty} (\isig,
\calo^n)$ for which $C_f$ satisfies condition $\eqq{1.2}$ is open in
$\Gamma _{L^\infty} (\isig, \calo^n)$. The map $\Phi: \calm^{(n)}_{L^\infty}
(V) \to \calpi^{(1,1)}(V)$, $\Phi(f)\deff \eta_{C_f}$, is continuous and
injective. The topology on the image $\Phi(\calm^{(n)}_{L^\infty}(V))$
coincides with the {\sl weak} topology of $\Gamma_{L^\infty} (\isig, \calo^n)$.
\endstate

\newgrph{Rem1.3}
\bf Remark. \rm For $C$ as in the lemma we shall call the corresponding
degree $n$ of $P$ the {\sl degree} of $C$.

\proof{1.2a} Since the group $\sfh^2(V,\zz)$ is trivial, any curve $C$ in $V$
admits a global determining function $F$. For $C$ satisfying condition 
$\eqq{1.2}$ the Weierstra\3 preparation theorem (see, \eg [GrHa]) insures
that $C$ is a zero-divisor of a uniquely defined Weierstra\3 polynomial
$P_f(z,w)=w^n + \sum _{i=1}^n f_i(z)\, w^{n-i}$, such that $F= h\cdot P_f$
for some {\sl invertible} $h\in \Gamma(V, \calo)$. We can
view $f \deff (f_1,\ldots,f_n)\in \Gamma(\isig, \calo^n)$ as a holomorphic map
from $\isig$ into the $n$-th symmetric power $\sym^n\Delta \subset \sym^n\cc
\cong\cc^n$. Hence the $f_i$ are necessarily uniformly bounded in $\isig$.
Moreover, for $g \in \Gamma_{L^\infty} (\isig, \calo^n)$ sufficiently close to
$f$, the curve $C_g$ defined by Weierstra\3 polynomial $P_g(z,w) \deff w^n + 
\sum_{i=1}^n g_i(z)\, w^{n-i}$, also satisfies condition $\eqq{1.2}$. Thus
the set $\calm^{(n)} _{L^\infty} (V)$ is open in $\Gamma_{L^\infty} (\isig,
\calo^n)$.

According to the Poincar\'e-Lelong formula (see [Ha] or [GrHa]), the map
$\Phi$ is given by formula
$$
f\in \Gamma_{L^\infty} (\isig, \calo^n) \longto
  \msmall{1\over \isl\pi}\,\d\dbar\,\log|P_f| \in \calpi^{(1,1)}(V).
$$
Thus $\Phi: \Gamma_{L^\infty} (\isig, \calo^n) \to \calpi^{(1,1)}(V)$ is
continuous.

Now let $f^{(\nu)} \in \calm^{(n)}_{L^\infty}(V) \subset \Gamma_{L^\infty}
(\isig, \calo^n)$ be a sequence. If $\{f^{(\nu)}\}$ converges weakly to $f\in
\calm^{(n)}_{L^\infty}(V)$, then $\{C_{ f^{(\nu)} } \}$ is bounded in $\calpi
^{(1,1)}(V)$ and any Cauchy subsequence of $\{C_{ f^{(\nu)} } \}$ must
converge to $C_f$. Vice versa, let $\{C_{ f^{(\nu)} } \}$ converge to a curve
$C$ satisfying condition $\eqq{1.2}$. Since $\calm^{(n)}_{L^\infty}(V)
\subset \Gamma_{L^\infty} (\isig, \calo^n)$ is bounded subset,
some subsequence of
$\{f^{(\nu)}\}$ converges weakly to $f\in \Gamma_{L^\infty} (\isig, \calo^n)$
such that $C_f =C$. Consequently, $f\in \calm^{(n)}_{L^\infty}(V)$.
\qed

\smallskip
\newgrph{Sm1.5}
We see that we have constructed a holomorphic family $\{C_f \}_{f\in
\calm^{(n)}_{L^\infty}(V)}$ of curves which are ``uniformly bounded" in $V$.
Later we shall show that this family possesses the universality property. To
generalise this result to other classes of boundary conditions on curves
we introduce the notion of a {\sl smoothness}.

For $0<r<R<\infty$ denote by $A_{r,R}$ an {\sl annulus} $\{\, z\in\cc\;:\;
r<|z|<R\,\}$. De\-note also $\Delta^-(r)\deff \{\, z\in \barr\cc \;:\; |z|>r
\,\}$. Recall that any Riemann surface $\isig$ with a  complex structure
which is homeomorphic to an annulus and which boundary $\d\isig$ consists of
two circles is in fact biholomorphic to some annulus $A_{r,R}$, $\isig \cong
A_{r,R}$. Recall also that any holomorphic function $f\in \Gamma(A_{r,R},
\calo)$ admits a unique decomposition into the sum $f=f^+ + f^-$ such that 
$f^+ \in \Gamma(\Delta(R), \calo)$ and $f^- \in \Gamma(\Delta^-(r), \calo)$ 
with $f^-(\infty)=0$. We call it the {\sl Laurent decomposition} of $f$.

\newdefinition{Smooth} \sl A smoothness class $S$ (\.\rm or simply a \sl
 smoothness) in $\Delta$ \rm is defined by fixing a subalgebra
$\Gamma_S(\Delta, \calo) \subset \Gamma_{L^\infty}(\Delta,\calo)$ which
satisfies the following conditions:

\nobreak
\slsi $\Gamma_S(\Delta, \calo)$ is a Banach algebra%
\footnote[1)]{
This means that $\norm{fg}_S\le c_S\norm{f}_S\norm{g}_S$,  where $c_S$
is a constant independent of $f,\,g\in \Gamma_S(\Delta, \calo)$ and 
possibly greater than 1. This can be always corrected by introducing a
new norm
$\norm{f}^*_S \deff \sup\left\{\,{\norm{fg}_S\over\norm{g}_S}\;:\; g\not=0\in 
\Gamma_S(\Delta, \calo)\,\right\}$ which is equivalent to $\norm\cdot$ 
and for which $\norm{fg}_S^*\le \norm{f}_S^*\norm{g}_S^*$.
}
with the norm $\norm{\cdot}_S$ and $\norm{f}_{L^\infty(\Delta)} \le C_S
\norm{f}_S$.

\slsii $\Gamma_S(\Delta, \calo)$ is invariant w.r.t.\ the action
of the group ${\bold U}(1)$ by rotations on $\Delta$.

\slsiii
If $f\in \Gamma_S(\Delta, \calo)$, $g\in\Gamma(A_{r,R},\calo)$ with some
$r<1<R$, and $fg= (fg)^+ + (fg)^-$ is the Laurent decomposition of the product,
then $(fg)^+ \in \Gamma_S(\Delta, \calo)$. Moreover $\norm{(fg)^+}_S \le
C(S,r,R)\cdot \norm{f}_S \cdot \norm{g}_{L^\infty(A_{r,R})}$ for some constant
$C(S,r,R)$.

\slsiv If $f\in \Gamma_{L^\infty}(A_{r,1}, \calo)$ has the Laurent
decomposition $f= f^+ + f^-$ with $f^+\in
\Gamma_S(\Delta, \calo)$ and a bounded invertible $g\deff {1/f} \in
\Gamma_{L^\infty}(A_{r,1}, \calo)$ with the Laurent decomposition $g=g^+ +
g^-$, then $g^+\in \Gamma_S( \Delta, \calo)$.

\smallskip
We say that $f\in \Gamma_S(\Delta, \calo)$ is {\sl $S$-smooth}.
Conditions \slsiii and \slsiv provides that $S$-smoothness depends
essentially only on the behavior of $f$ at the boundaryof $\Delta$.
Obvious examples are $C^k$-Lipschitz-H\"older smoothness (up to the
boundary) $S=C^{k,\alpha}(\barr\Delta)$ with $k\in\nn$ and $0\le\alpha\le1$,
Sobolev smoothness $S=L^{k,p}(\Delta)$ with $k\ge1$, $1\le p\le \infty$ and
$k\,p>2$, and also Sobolev smoothness $S=L^{k,p}(S^1)$ on boundary $S^1\deff
\d\Delta$ with $k\ge1$ and $1\le p\le \infty$. The later means that the trace
$f\ogran_{S^1}$ of $f\in \Gamma(\Delta, \calo)$ on $S^1=\d \Delta$ 
is well defined and belongs to the corresponding class.

\newdefinition{1.6} Let $\isig$ be a smooth complex curve whose boundary
$\d\isig$ consists of finitely many components $\gamma_i$, $i=1,\ldots,n$,
each of which is homeomorphic to a circle $S^1$. {\sl The smoothness $S$ at
$\d\isig$ in $\isig$} is defined by fixing of a smoothness classes $S_i$ in
$\Delta$ and annuli $A_i\subset \isig$ such that one of the components of
boundary $\d A_i$ coincides with $\gamma_i$ and the other one lies in the
interior of $\isig$. For every $i=1,\ldots,n$ this induces a biholomorphic map
$\phi_i: A_i \to A_{r_i,1}$ which extends continuously
up to the boundary $\d A_i$ and
maps $\gamma_i$ onto $\d\Delta$. We say that $f\in \Gamma(\isig, \calo)$
is {\sl $S$-smooth in $\isig$ at $\d\isig$}, $f\in \Gamma_S(\isig,\calo)$,
\iff for every $i=1, {\ldots},n$ the Laurent decomposition $\phi_{i\,*}f=
(\phi_{i\,*} f)^+ + (\phi_{i\,*}f)^-$ yields $(\phi_{i\,*}f)^+ \in
\Gamma_{S_i}(\Delta, \calo)$.

\newlemma{?1} Let $\isig$, $\gamma_i$, $\phi_i: A \to A_{r_i,1}$, and $S_i$
be as above. Then $\Gamma_S(\isig,\calo)$ is a Banach algebra with respect to
the norm $\norm{f}_S\deff \sum_i \norm{(\phi_{i\,*}f)^+}_{S_i}$. Moreover,
$$
\norm{f}_{L^\infty(\isig)} \le C_S\cdot \norm{f}_S
\neqq{N0}
$$
with a constant $C_S$ independent of $f\in \Gamma_S(\isig,\calo)$. The subset
$$
\Gamma_S(\isig,\calo)^\times \deff \{\,f\in \Gamma_S(\isig,\calo) \;:\;
f\inv\in \Gamma_{L^\infty}(\isig,\calo) \,\}
$$
is open in $\Gamma_S(\isig,\calo)$ and the map
$$
F: \Gamma_S(\isig,\calo)^\times \to \Gamma_{L^\infty}(\isig,\calo)
\qquad F(f)\deff f\inv
$$
is $\Gamma_S(\isig,\calo)$-valued and holomorphic.
\endstate

\proof{?1p} Let $f\in \Gamma_S(\Delta,\calo)$. Set $f_i \deff \phi_{i\,*}f$
and let $f_i = f_i^+ + f_i^-$ be the corresponding Laurent decompositions.
Due to \slsi of the definition of smoothness, $f$ is uniformly bounded in
$\isig$ and takes its supremum on one of the boundary component $\gamma_j$.
The later means that there exists a sequence $x_\nu\in \isig$ such that
$\lim x_\nu \in \gamma_j$ and $\lim |f(x_\nu)| = \norm{f}_{L^\infty(\isig)}$.
In particular, $\norm{f}_{L^\infty(\isig)} =\norm{f_j}_{L^\infty(A_j)}$.

For this $A_j$ we obviously have
$$
\norm{f_j^-}_{L^\infty(A_j)} \le \norm{f_j}_{L^\infty(A_j)} +
\norm{f_j^+}_{L^\infty(A_j)}.
\neqq{N1}
$$

On the other hand, $\norm{f_j^-}_{L^\infty(\gamma_j)}\le \delta \cdot
\norm{f_j^-}_{L^\infty(A_j)}$ with $\delta\deff\max\{r_i\}<1$. This is due to
$f_j^-(\infty)=0$ and the Schwarz inequality. Consequently,
$$
\norm{f_j}_{L^\infty(A_j)} \le \delta\cdot \norm{f_j^-}_{L^\infty(A_j)} +
\norm{f_j^+}_{L^\infty(A_j)}.
\neqq{N2}
$$
Comparing $\eqq{N1}$ and $\eqq{N2}$ we see that
$$
\norm{f_j^-}_{L^\infty(A_j)} \le {2\over1-\delta}\cdot
\norm{f_j^+}_{L^\infty(A_j)} \quad\text{and}\quad
\norm{f_j}_{L^\infty(A_j)} \le {3-\delta\over1-\delta}
\cdot \norm{f_j^+}_{L^\infty(A_j)} .
\neqq{N3}
$$

Since $\delta<1$ is independent of $f$, from \slsi we obtain the estimate
$\eqq{N0}$. Consequently, every$\norm{\cdot}_S$-Cauchy sequence converges
to some element in $\Gamma_S(\isig,\calo)$. Moreover, the map
$$
\Phi: \Gamma_S(\isig,\calo) \longto
\prod_{i=1}^n \Gamma_{S_i}(\Delta,\calo),  \qquad
\Phi(f)\deff\bigl((\phi_{i\,*}f)^+\bigr),
\neqq{N4}
$$
is a closed imbedding.

\smallskip
Take another $g\in \Gamma_S(\isig,\calo)$ and set $g_i \deff \phi_{i\,*}g$
with the corresponding Laurent decompositions $g_i = g_i^+ + g_i^-$. Then
$$
(\phi_{i\,*}(f\cdot g))^+ = f_i^+ g_i^+ + (f_i^- g_i^+)^+ +
(f_i^+  g_i^-)^+ .
$$
For any $i=1,\ldots,n$ we obviously have
$$
\norm{f_i^-}_{L^\infty(\Delta^-(r_i))}
\le c(r_i) \cdot \norm{f_i}_{L^\infty(A_{r_i,1})}
\le c_S \norm{f}_S
$$
and the same estimates for $g$. Due to \slsiii, we obtain
$$
\norm{f\cdot g}_S \le c_S'\cdot \norm{f}_S \cdot \norm{g}_S.
$$

\smallskip
Condition \slsiv implies that $\Gamma_S(\isig,\calo)^\times$ consists of
those $f\in \Gamma_S(\isig,\calo)$ which are invertible in $\Gamma_S(\isig,
\calo)$. Thus the last statement of the theorem is a standard fact of the
theory of commutative Banach algebras, see \eg [Ga].
\qed

\smallskip
\newdefinition{FisSsmooth} Let $\isig$ and $S$ be as above and $V=\isig \times
\Delta$. A function $F(z,w)\in \Gamma(V,\calo)$ is called {\sl $S$-smooth},
$F\in \Gamma_S(V,\calo)$, \iff $F(z, a)\in \Gamma_S(\isig,\calo)$ for every
$a\in \Delta$ and the induced map $F: \Delta \to \Gamma_S( \isig, \calo)$ is
holomorphic and bounded. A curve $C\subset V$ is called {\sl $S$-smooth} \iff
$C$ satisfies  condition $\eqq{1.2}$ and $C = \div(F)$ for some $F\in
\Gamma_S(V, \calo)$, \ie $C$ is defined by such an $F$. Let $\calf_S^{(n)}$
denote the set of those $F\in \Gamma_S(V, \calo)$ such that $\inf\{\,|F(z,a)|
\;:\; z\in\isig, a\in A_{r,1}\,\} >0$ for some $r<1$ and for which the curve
$C_F\deff\div(F)$ has degree $n$.

\newlemma{?2} A set $\Gamma_S(V,\calo)$ of $S$-smooth functions is a Banach
space with a norm
$$
\norm{F}_S \deff \sup\{\,\norm{F(z,a)}_S \;:\; a\in\Delta\,\}.
$$
Every $S$-smooth curve $C\subset V$ of degree $n$ is represented by a unique
Weierstra\3 polynomial $P= P_f(z,w)= w^n + \sum_{i=1}^n f_i(z) w^{n-i}$ with
$f= (f_1, \ldots, f_n) \in \Gamma_S(\isig, \calo^n)$. The set $\calf_S^{(n)}$ 
is open in $\Gamma_S(V, \calo)$ and the map $\Psi: \calf_S^{(n)} \to 
\Gamma_S(\isig, \calo^n)$, $\Psi:F\mapsto f$, is holomorphic.
\endstate

\proof{?2} The part of the lemma concerning $\Gamma_S(V,\calo)$ is obvious.
Let $F$ lie in $\calf_S^{(n)} \subset \Gamma_S(V, \calo)$ and $C_F\subset \isig
\times\Delta$ be the corresponding curve. For $k\in\nn$ and $z\in\isig$ we
set
$$
f_k(z) \deff \msmall{1\over 2\pi\isl} \int_{|w|=r} \msmall{w^k\over F(z,w)}
\msmall{\d F\over \d w}(z,w) \,dw ,
$$
where $r<1$ is chosen sufficiently close to 1. Then $f_0$ is constant and
equals the degree $n$ of $C$, whereas $f_i$, $i=1, \ldots,n$, are the
coefficient of the Weierstra\3 poly\-nomial $P_f$ of $C$. 
Since the operations
in the definition of $f_k$---taking inverse, differentiating,
integrating---are holomorphic, $f_k$ depends holomorphically on $F\in
\Gamma_S(V, \calo^n)$.
\qed

\smallskip
\newgrph{?2a} Let $\calm^{(n)}_S(V)$ be the image $\Psi(\calf_S^{(n)})$. One
can regard an embedding $\calm^{(n)}_S(V) \subset \Gamma_S(\isig, \calo^n)$ 
as a local chart of a moduli space of curves on a complex surface with an
appropriate smoothness condition. To be able to patch such local
models together we need an invariant description of $\calm^{(n)}_S(V)$.

Take $F\in \calf_S^{(n)}$, set $C\deff\div(F)$, $f\deff\Psi(F)$, 
and consider the tangent map
$$
d\Psi_F: T_F\calf_S^{(n)}=\Gamma_S(V,\calo) \longto
T_C \calm^{(n)}_S(V)=\Gamma_S(\isig,\calo^n).
$$

\newlemma{?3} Let $\pi:|C|\to \isig$ be the natural projection and $\pi_*
\calo_C$ a push-forward with respect to $\pi$. The map $d\Psi_F$ induces on
$\calo_\isig^n$ a structure of a $\pi_*\calo_C $-module which is independent
of the choice of $F$ with $\div(F)=C$. With respect to this structure there
exists a $\pi_*\calo_C $-isomorphism $\theta_C$ between $\calo_\isig^n$ and
the push-forward $\pi_*\caln_C$ of the normal sheaf $\caln_C = \Hom_{\calo_V}
\vph (\cali_C/\cali_C^2, \calo_C)$.

The isomorphism $\theta_C$ admits the following characterization: Let
$F_\lambda(z,w)$, $\lambda\in \Delta(\rho)$, be a holomorphic family of
functions in $\calf_S^{(n)}$ such that $C=\div(F_0)$. Furthermore, let $\phi
\deff \theta_C(d\Psi_{F_0}(F'_0))\in \Gamma(C, \caln_C)$ where $F'_0\deff 
{\d F_\lambda\over \d \lambda}\ogran_{\lambda=0}$. Then
$$
\phi: [F_0]_{\cali_C^2} \in \Gamma(C, \cali_C/\cali_C^2) \mapsto
[F'_0]_{\cali_C}\in \Gamma(C, \calo_V/\cali_C).
\neqq{theta}
$$
\endstate

\proof{?3p} For convenience we slightly modify the definition of the map
$\Psi$. For $F\in \calf_S^{(n)} \subset \Gamma_S(V,\calo)$ and $f=(f_1,\ldots,
f_n) =\Psi(F) \in \Gamma_S(\isig,\calo^n)$ we set $\Psi^{W\!P}(F)\deff w^n +
\sum_{i=1}^n f_i(z) w^{n-i}\in \Gamma_S(V,\calo)$, where  $W\!P$
stands for ``Weierstra\3 polynomial".

Set $P\deff \Psi^{WP}(F)$ and $g\deff P\inv F$. 
Then $g$ is holomorphic and bounded in $V$, $g(z,a)\in
\Gamma_S(\isig, \calo)$ for $a\in A_{r,1}$, and the induced map $g: A_{r,1}
\to \Gamma_S(\isig, \calo)$ is holomorphic. Considering the Cauchy
representation for $g$,
$$
g(z,w)=\msmall{1\over 2\pi\isl} \int_{|\zeta|=r}
\msmall{g(z,\zeta)\over w-\zeta}\,d\zeta,
$$
we see that $g\in \Gamma_S(V, \calo)$. By the same argumentation $g\inv \in
\Gamma_S(V, \calo)$. Moreover, the map $F\in \calf_S^{(n)} \longmapsto
F/\Psi^{W\!P}(F)\in \Gamma_S(V, \calo^n)$ is holomorphic.

Now let $F_\lambda(z,w)$, $\lambda\in \Delta(\rho)$, be a holomorphic family
of functions in $\Gamma_S(V, \calo)$ with $F_0=F$. Put $P_\lambda \deff
\Psi^{W\!P}(F_\lambda)$ and $g_\lambda \deff F_\lambda /P_\lambda$. Then
$P_\lambda$ and $g_\lambda$ depend holomorphically in $\lambda$.
Differentiating the identity $F_\lambda = g_\lambda \cdot P_\lambda$ with
respect to $\lambda$ in $\lambda=0$, we obtain
$$
F'_0 = g'_0 \cdot P_0 + g_0 \cdot P'_0
$$
with $g_0=g$ and $P_0=P$. Thus the tangent map
$$
d\Psi^{W\!P}_F : T_F\calf_S^{(n)} =\Gamma_S(V,\calo) \to \Gamma_S(V,\calo)
$$
is given by a formula $d\Psi^{W\!P}_F(F')= \left\{ g\inv F'\over P \right\}$,
the Weierstra\3 remainder of a division of $g\inv F'$ by $P$. It is a unique
polynomial $R_{f'}= \sum_{i=1}^n f'_i(z) w^{n-i}$ of degree $<n$ such that 
$R_{f'} \equiv g\inv F'(\mod F)$. This implies that $d\Psi_F$ yields an 
isomorphism of Banach spaces
$$
\eqalign{
\psi_F:
&  \Gamma_S(V, \calo) / \bigl(F\cdot \Gamma_S(V, \calo) \bigr)
\  \cong\ \Gamma_S(\isig, \calo^n),
\cr
\psi_F:
&
[H]_F \longmapsto h=(h_1,\ldots,h_n)
\text{ with } g\inv H\equiv R_h (\mod F).
}
$$
Due to its definition, $\psi_F$ is essentially local and induces the
isomorphism of $\calo_\isig$-modules
$$
\psi_F : \pi_*\bigl(\calo_V/ (F(z,w) \cdot \calo_V) \bigr)
\cong
\calo_\isig[w]/(P(z,w) \cdot \calo_\isig[w] )
\cong \calo_\isig^n.
$$
Since $\calo_V/ P(z,w) \cdot \calo_V = \calo_C$, $\psi_F$ defines on
$\calo_\isig^n$ a structure of a free $\pi_* \calo_C $-module of rank 1. If
$\wt F\in \Gamma_S(V,\calo)$ is another function with $\div(\wt F)=C$, then
$\wt F= h\cdot F$ for some invertible $h\in\Gamma_S(V,\calo)$. Then, by the
definition, $\psi_{\wt F}(H) = \psi_F( h\inv H)$. Consequently, the induced
structure of the $\pi_* \calo_C$-module on $\calo^n_\isig$ is independent
of the choice of $F$.

\medskip
Using $\psi_F$ we define a $\pi_* \calo_C $-homomorphism $\theta_F:
\calo_\isig^n \to \pi_* \caln_C$. For a local section $f=(f_1, \ldots,
f_n)$ of $\calo_\isig^n$ over an open set $\iom\subset \isig$ we take a
holomorphic function $H\in \Gamma(\iom\times \Delta, \calo)$ such that $H
\equiv g\cdot R_f (\mod F)$ in $\iom\times \Delta$, where $g= P\inv F$ is
as above. Since $C=\div(F)$, the sheaf $\cali_C/ \cali_C^2$ is free
$\calo_C$-module of rank 1 with generator $[F]_{\cali_C^2}$. Using
$\caln_C  = \hom_{\calo_V}(\cali_C/\cali_C^2, \calo_V/\calo_C )$,
we define 
$$
\theta_F(f)  \in \Gamma(\iom, \pi_*\caln_C)
  = \Gamma\bigl( (\iom\times\Delta)\cap C, \caln_C\bigr), \qquad
\theta_F(f)  :\,[F]_{\cali_C^2} \longmapsto [H]_{\cali_C}.
$$

If $\wt H\in \Gamma(\iom\times \Delta, \calo)$ is another holomorphic
function with $\wt H\equiv g\cdot R_f (\mod F)$, then $ [\wt H]_{\cali_C} =
[H]_{\cali_C}$. This shows that the definition of $\theta_F(f)$ is independent
of the choice of $H$.

Similarly, if $\wt F\in \Gamma_S(V,\calo)$ is another defining function for
$C$, $C=\div(\wt F)$, then $\wt F =h\cdot F$ with $h \in \Gamma_S(V,\calo)$
invertible. In this case $P\inv \wt F= h\cdot g$ and hence $\theta_{\wt
F}(f): [\wt F]_{\cali_C^2} \mapsto [h\cdot H]_{\cali_C}$. This means that
$\theta_{\wt F}(f) = \theta_F(f)$ as sections of $\caln_C$. 
Thus the definition of $\theta_C \deff \theta_F$ is independent of the choice 
of defining function $F$ for the curve $C$.

\smallskip
Now let $F_\lambda(z,w) \in \Gamma_S(V,\calo)$, $\lambda\in \Delta(\rho)$, be
a holomorphic family of functions such that $F_0=F$ and $F'_0=H$. The relation
$\eqq{theta}$ for $\phi\deff \theta_C(d\Psi_F(H))$ follows immediately from
the construction of $\Psi_F$ and $\theta_C$.
\qed

\newdefinition{GsN}
Note that the constructed isomorphism $\theta_C$ induces the bijection
$\theta_C: \Gamma(C, \caln_C) \cong \Gamma(\isig, \calo_\isig^n)$. A
section $\phi$ of $\caln_C$ is called {\sl $S$-smooth}, $\phi \in \Gamma_S(C,
\caln_C)$, \iff $\theta(\phi) \in \Gamma_S(\isig, \calo_\isig^n)$.

\newgrph{GsN1}
As we have already noted, the map $\Psi: \calm_S^{(n)}(V) \to \Gamma_S(\isig,
\calo_\isig^n)$ is a chart in a moduli space of curves in $V$ with an
appropriate smoothness condition at boundary. However, the map $\Psi$
depends on the choice of local (holomorphic) coordinates $(z,w)$ in $V$. In
particular, if in the construction of $\Psi$ we replace $w$ by some other
coordinate function $\ti w = \ti w(w)$, then the map $\Psi$ as well as
a $\pi_* \calo_C $-module structure on $\calo_\isig^n$ will change.
However, relation $\eqq{theta}$ remains valid, since it is independent
of the choice of coordinates in $V$. This leads us to

\newcorollary{TcM} The tangent space $T_C\calm_S^{(n)}(V)$ is canonically
isomorphic to $\Gamma_S(C, \caln_C)$.
\qed
\endstate

\bigskip\bigskip
\newsection{S2}{Curves in the ``distorted cylinder"}

\nobreak
\newgrph{2.1}
To be able to patch local descriptions, we consider the following special
situation.
In $\cc^2$ with  standard coordinates $(z,w)$ we consider
an annulus $A_{r,R}\deff \{\, (z,w) \;:\; w=0,\; r<|z|<R \,\}$. We assume that
in some neigh\-bour\-hood $U$ of the closure $\barr A_{r,R}$ we are given
two holomorphic functions $z_1$ and $z_2$ which {\sl coincide with $z$ along
$A_{r,R}$}. Without loss of generality we may also assume that the both pairs
$(z_1, w)$ and $(z_2, w)$ are coordinates in $U$ so that we can express
$z_1=z_1(z_2, w)$ and $z_2=z_2(z_1, w)$.

For $\rho>0$ let
$$
W_{r,R,\rho} \deff
  \{\,x\in U\;:\; |z_1(x)|>r, |z_2(x)|<R, |w(x)|<\rho \,\}.
$$
We shall always suppose that $\rho>0$ is chosen sufficiently small 
such that $W_{r,R,\rho} \Subset U$ and that the sets
$$
\eqalign{
\d_-W_{r,R,\rho}
& \;\deff\;
\{\,x\in U\;:\; |z_1(x)|=r,  |w(x)|\le\rho \,\}
\cr
\d_+W_{r,R,\rho}
& \;\deff\;
\{\,x\in U\;:\; |z_2(x)|=R, |w(x)|\le\rho \,\}
}
$$
are disjoint. One can regard the set $W_{r,R,\rho}$ as a distorted cylinder
with the non-parallel lower side $\d_-W_{r,R,\rho}$ and upper side
$\d_+W_{r,R,\rho}$. Note also that there exist real numbers $r<r'<R'<R$
such that both  sets
$$
\eqalign{
V^-_{r,r'}
& \;\deff\;
 \{\,x\in U\;:\; r<|z_1(x)|<r',  |w(x)|<\rho
\,\}
\cr
V^+_{R',R}
& \;\deff\;
\{\,x\in U\;:\; R'<|z_2(x)|< R,  |w(x)|<\rho \,\}
}
$$
are products of an annulus and a disk. This allows us to make the following

\newdefinition{2.2} {\sl A smoothness} $S$ in $W_{r,R,\rho}$ is defined by
fixing smoothness classes $S^-$ and $S^+$ in $\Delta$. A holomorphic function
$F$ in $W_{r,R,\rho}$ is {\sl $S$-smooth}, $F\in \Gamma_S(W_{r,R,\rho},
\calo)$, \iff $F \ogran_{V^-_{r,r'}} \in \Gamma_{S^-}(V^-_{r,r'}, \calo)$ and
$F \ogran_{V^+_{R',R}} \in \Gamma_{S^+}(V^+_{R',R},\calo)$. Note that $S$
defines also a smoothness in $A_{r,R}$: It is suffucient to fix annuli 
$A^-=\{\, r<|z|<r' \,\}$, $A^+=\{\, R'<|z|<R \,\}$ and smoothness classes 
$S^+$ and $S^-$. We shall also denote this smoothness by $S$. Thus  we obtain 
a continuous projection map $F \in \Gamma_S(W_{r,R,\rho}, \calo) \mapsto 
F\ogran_{A_{r,R}} \in \Gamma_S(A_{r,R}, \calo)$.

A curve $C\subset W_{r,R,\rho}$ is {\sl $S$-smooth} \iff $C \subset
W_{r,R,\rho'}$ for some $\rho'< \rho$ and $C= \div(F)$ for some $F\in \Gamma_S
(W_{r,R,\rho}, \calo)$. The {\sl degree of an $S$-smooth curve $C\subset
W_{r,R,\rho}$} is an integer $\deg C\deff \int_\gamma d\log F$, where $F\in
\Gamma_S( W_{r,R,\rho}, \calo)$ is any function defining $C$ and $\gamma$ is
a simple smooth loop (\ie a closed real curve) in $\{\, x\in W_{r,R,\rho}
\;:\; |w(x)|=\rho \,\}$ with $\int_\gamma d\log w=1$. It is obvious that
$\deg C$ is a positive interger independent of the choice of $F$ and $\gamma$.
The set of $S$-smooth curves of degree $n$ in $W_{r,R,\rho}$ will be denoted
by $\calm^{(n)}_S(W_{r,R,\rho})$. Note that for every $S$-smooth curve $C
\subset W_{r,R,\rho}$ the Weierstra\3 polynomials $P^+$ and $P^-$ of $C \cap
V^-_{r,r'}$ and $C \cap V^+_{R',R}$ are uniquely defined. This yields an
injective map
$$
\kappa^{(n)}: \calm^{(n)}_S(W_{r,R,\rho}) \to
\Gamma_{S^-}(A_{r,r'}, \calo^n) \times
\Gamma_{S^+}(A_{R',R}, \calo^n).
\neqq{kappa}
$$

A family $\{C_y\}_{y\in Y}$ of $S$-smooth curves of a degree $n$ in $W_{r,R,
\rho}$ parametrised by a topological space $Y$ is called {\sl continuous} 
\iff induced map
$$
\Psi : Y \to \Gamma_{S^-}(A_{r,r'}, \calo^n) \times
\Gamma_{S^+}(A_{R',R}, \calo^n), \qquad \Psi(y)\deff \kappa^{(n)}(C_y),
$$
is continuous. 

We say that 
$\{C_y\}_{y\in Y}$  is a {\sl holomorphic family of $S$-smooth curves}, if 
it is a continuous family of $S$-smooth curves, $Y$ has a structure of a 
Banach analytic space, and both restricted families $C_y \cap V^-_{r,r'}$ and 
$C_y \cap V^+_{R',R}$ are given by holomorphic morphisms $\bfpsi^-_Y: Y \to 
\Gamma_{S^-}(V^-_{r,r'}, \calo)$ and $\bfpsi^+_Y: Y \to \Gamma_{S^+}
(V^+_{R',R}, \calo)$. Note that $\bfpsi^\pm_Y$ induce local determining 
functions $F^\pm(z,w,y) \deff \psi^\pm_Y(y)(z,w)$ on $V^-_{r,r'}\times Y$ and 
$V^+_{R',R}\times Y$ respectively.

\smallskip
To generalise the results of {\sl Section 1\/} for curves in $W_{r,R,\rho}$
one must find an appropriate analog of a Weierstra\3 polynomial for
$W_{r,R,\rho}$.

\newdefinition{dWP} Let the components of $f=(f_1,\ldots,f_n)\in
\Gamma_S(A_{r,R}, \calo^n)$ have a Laurent decomposition $f_i(z)= f_i^+(z) +
f_i^-(z)$. A {\sl (distorted) Weierstra\3 polynomial $\wt P_f(z_1,w)$ in
$W_{r,R,\rho}$ of degree $n$ with coefficients $(f_1,\ldots,f_n)$} is defined
as
$$
\wt P_f(z_1,z_2,w) \deff
w^n + \sum_{i=1}^n \bigl(f^-_i(z_1)+f^+_i(z_2) \bigr)\,w^{n-i}.
\neqq{2.5}
$$

\newgrph{2.0}
One can expect that there is one-to-one correspondence between $S$-smooth
curves in $W_{r,R,\rho}$ of degree $n$ and distorted Weierstra\3 polynomials
$\wt P_f(z_1,z_2,w)$ of the same degree $n$ with $S$-smooth coefficients.
Since a difference of $z_1$ and $z_2$ introduces a ``non-linearity''
one can hope to obtain the corresponding relationship
only in some neighbourhood of a trivial case $\wt P_0=w^n$ and $C_0=n\cdot
A_{r,R}$, when the coefficients $f=(f_1,\ldots,f_n)$ of $\wt P_f$ are
sufficiently small with respect to the norm in $\Gamma_S(A_{r,R}, \calo^n)$.
The corresponding condition on a curve is that $\norm{\kappa^{(n)}(C)}_S$
(with $\kappa^n$ from \eqq{kappa})
should be small. Here $\norm{\cdot}_S$ denotes the norm in
$\Gamma_{S^-}(A_{r,r'}, \calo^n) \oplus \Gamma_{S^+}(A_{R',R}, \calo^n)$.

\newlemma{2.1}Let $z_1$, $z_2$, $w$, $r<r'<R'<R$, $\rho$, and $S$ have the
same meaning as above. There exists $\epsi>0$ such that every $S$-smooth curve
$C$ of degree $n$ in $W_{r,R,\rho}$, satisfying
$$
\norm{\kappa^{(n)}(C)}_S \le \epsi
\neqq{W-epsi}
$$
is a zero divisor of a uniquely defined distorted Weierstra\3 polynomial $\wt
P_f(z_1,z_2,w)$.

If $Y$ is a topological (resp.\ Banach analytic) space and $\{C_y\}_{y\in Y}$ 
is a continuous (resp.\ holomorphic) family of curves satisfying $\eqq{W-epsi}
$, then the induced map $\psi_Y :Y \to \Gamma_S(A_{r,R}, \calo^n)$ with
$\div(\wt P_{\psi_Y(y)}) =C_y$ is continuous (resp.\ holomorphic).
\endstate

\proof{2.1p} Without loss of generality we may assume that $\epsi>0$ is chosen 
sufficiently small so that $C \subset W_{r,R,\rho'}$ for any given $\rho'<
\rho$ also small enough. Furthermore, we may also assume that
$$
\matrix
r\;<|z_1(z_2,w)|<R',&
    \quad\hbox{for any $(z_2,w)$ with $|z_2|=r'\;$ and $|w|\le\rho'$,}
\cr
r'<|z_2(z_1,w)|<R,\;&
    \quad\hbox{for any $(z_1,w)$ with $|z_1|=R'$ and $|w|\le\rho'$.}
\endmatrix
\neqq{r'}
$$
Then $W_{r,R,\rho'}$ is a union of $V_1 \deff V^-_{r,R'} = \{\, r<|z_1|<R', \;
|w|< \rho' \,\}$ and $V_2 \deff V^+_{r',R} = \{\, r'<|z_2| <R,\; |w|< \rho'
\,\}$. Since $V_i= A_i\times \Delta(\rho')$ with $A_1\deff \{\, r<|z|<R' \,\}$
and $A_2\deff\{\, r'<|z|<R \,\}$, we can apply the result of \lemma{1.2}.

Let $B$ be a sufficiently small ball in $\Gamma_S(A_{r,R}, \calo^n)$, $f=
(f_1, \ldots, f_1)\in B$, and $\wt P_f$ a corresponding distorted Weierstra\3
polynomial. Then the zero divisor $C_f\deff \div (\wt P_f)$ is an $S$-smooth
curve of degree $n$ lying in $W_{r,R,\rho}$. Moreover, both curves $C_f \cap
V_i$, $i=1,2$, are $S_i$-smooth and of degree $n$ in $V_1$ and $V_2$
respectively. Here we set $S_1 \deff S^-$ and $S_2 \deff S^+$. Thus there 
exist uniquely defined Weierstra\3 polynomials $P_1$
and $P_2$ in $V_1$ and $V_2$ respectively such that $C_f \cap V_i =\div(P_i)$.
This defines the maps $\phi_i : B \to \Gamma_{S_i}(A_i, \calo^n)$, $i=1,2$.

Formula $\eqq{theta}$ provides that the derivation of $\phi_i$ at
$f\equiv 0\in B$ is simply the restriction map $\Gamma_S(A_{r,R}, \calo^n)
\to \Gamma_{S_i} (A_i, \calo^n)$. Set $Y \deff \Gamma_{S_1}(A_1, \calo^n)
\oplus \Gamma_{S_2} (A_2, \calo^n)$ and $\phi=(\phi_1,\phi_2): B\to Y$ so
that $\phi(f)=\kappa^{(n)}(\div(\wt P_f))$. Since the differential $d\phi(0)$
of $\phi$ at $f\equiv 0\in B$ consists of the pair of restrictions, $d\phi(0)$ 
is an injection with a closed image. This implies the injectivity 
of $\phi$ in some smaller ball $B(0,\epsi) \subset \Gamma_S(A_{r,R}, \calo^n)$.

\smallskip
We state our conclusion in the following way:
{\sl There exists an $\epsi>0$ such that two distorted Weierstra\3 polynomials
$\wt P_f$ and $\wt P_g$ in $W_{r,R,\rho}$ of given degree $n$ with
$\norm{f}_S \le \epsi$ and $\norm{g}_S \le \epsi$ coincide provided
they define the same curve $C\subset W_{r,R,\rho}$}.

\medskip\newgrph{2.1a}
Now let $C\subset W_{r,R,\rho}$ be a curve which satisfies the hypotheses
of the lemma. In particular,
$C\cap V_1$ is $S^-$-smooth and $C\cap V_1 = \div(P)$
for a uniquely defined Weierstra\3 polynomial $P=w^n + \sum_{i=1}^n
g_i(z_1)\, w^{n-i}$. Further, from $C\subset W_{r,R,\rho'}$ we obtain
$\norm{g_k}_{L^\infty (A_{r,R'})} \le c \cdot \rho'$ where  the constant $c$
is independent of the curve $C$. This yields
$$
\norm{g_k}_{\Gamma_{S^-}(A_{r,r'} , \calo)} \le c' \cdot \rho'.
\neqq{P0-est}
$$

Consider the restriction of $P$ to the set
$$
W_{r,r',\rho} = \{\,x\in U\;:\; |z_1(x)|>r, |z_2(x)|<r', |w(x)|<\rho \,\}.
$$
Note that every $S$-smooth function $F$ in $W_{r,r',\rho}$ is uniquely
represented in the form
$$
F= \sum_{i=1}^n (f^+_i(z_1) + f^-_i(z_2))\, w^{n-i} + w^n (1+Q)
\neqq{F=P+Q}
$$
with $Q\in \Gamma_S(W_{r,r',\rho}, \calo)$ and $f=(f_1, \ldots, f_n) \in
\Gamma_S(A_{r,r'}, \calo^n)$. Here $f_i(z)= f^+_i(z) + f^-_i(z)$ denotes the
Laurent decomposition of the components of $f(z)$.
The corresponding $f_i$ are obtained inductively by the formula
$$
f_{n-k}(z) \deff \msmall{F- \sum_{i=n-k+1}^n (f^+_i(z_1) + f^-_i(z_2))\,
w^{n-i} \over w^k}\ogran\vph_{A_{r,r'}} \qquad k=0, \ldots, n-1,
\neqq{Fk}
$$
so that
$$
1 + Q \deff \msmall{F- \sum_{i=1}^n (f^+_i(z_1) + f^-_i(z_2))\,
w^{n-i} \over w^n}\cdot
\neqq{Q}
$$

Let us denote $P$ by $\wt P_0$. Define inductively $\wt P_{k+1}\deff
(1-Q_k)\inv \wt P_k$, where $Q_k$ is determined by the relation
$$
\wt P_k= \sum_{i=1}^n (f^+_{k,i}(z_1) + f^-_{k,i}(z_2))\, 
w^{n-i} + w^n (1+Q_k),
\neqq{Pk}
$$
and $f_{k,i}(z) = f^+_{k,i}(z) + f^-_{k,i}(z)$ is the Laurent
decomposition. We shall represent $\wt P_k$ in the form $\wt P_k= w^n+ R_k+
w^n Q_k$ with $R_k(z_1,z_2,w)\deff \sum_{i=1}^n (f^+_{k,i}(z_1) + f^-_{k,i}
(z_2)) \, w^{n-i}$. The estimates $\eqq{P0-est}$ on the coefficients $g_i$ of
$\wt P_0$ and the recursive formulas for $f_{0,i}(z)$ and $Q_0$ provide the
estimate
$$
\norm{R_0} + \norm{Q_0}
\le c'' \cdot \rho'\!,
$$
where $\norm{\cdot}$ denotes the norm in $\Gamma_S(W_{r,r',\rho}, \calo)$
and the constant $c''$ independent of the choice of a curve $C$. In the same
way one obtains the estimate
$$
\norm{R_{k+1}} \le (1 + c'''\,\norm{Q_k})\,\norm{R_k}
\quad\hbox{and} \quad
\norm{Q_{k+1}} \le c'''\,(\norm{Q_k} + \norm{R_k})
\,\norm{Q_k},
$$
where the constant $c'''$ is independent of $C$ and $\rho'$. Since $Q_0$ and
$R_0$ are small enough, the iteration converges to $\wt P_\infty=
w^n+R_\infty$ which is of the desired form.

\smallskip
To show the existence of a distorted Weierstra\ss-type polynomial $\wt P$ in
the entire set $W_{r,R,\rho}$, we additionally fix real numbers $r''$ and $R''$
with the property $r'<r''<R''<R'$. Then there exists a $\rho''>0$ such that
$$
\matrix
r'\;<|z_2(z_1,w)|<R'',&
    \quad\hbox{for any $(z_1,w)$ with $|z_1|=r''\;\,$ and $|w|\le\rho''$,}
\cr
r''<|z_1(z_2,w)|<R',\;&
    \quad\hbox{for any $(z_2,w)$ with $|z_2|=R''$ and $|w|\le\rho''$.}
\endmatrix
\neqq{r''}
$$
We may assume that $\epsi>0$ was fixed so small that every curve $C$
satisfying the hypotheses of the lemma lies in $W_{r,R,\rho''}$.

Let $C$ be such a curve. The above procedure allows us to construct the
corresponding distorted Weierstra\ss-type polynomials $\wt P^-$ in the set
$W_{r,R'',\rho''}$ and $\wt P^+$ in the set $W_{r'',R,\rho''}$. Due to
condition $\eqq{r''}$ the intersection $W_{r,R'',\rho''} \cap W_{r'',R,\rho''}$
is also a distorted cylinder $W_{r'',R'',\rho''}$. Thus $\wt P^-$ and
$\wt P^+$ coincide and define the desired distorted Weierstra\ss-type
polynomial $\wt P$ in the whole set $W_{r,R,\rho}$.

\newgrph{univ}
Let $Y$ be a topological space and $\{C_y\}_{y\in Y}$ a continuous family 
of curves satisfying condition $\eqq{W-epsi}$. Furthermore, let $\psi_Y :Y \to 
\Gamma_S(A_{r,R}, \calo^n)$ be an induced map. The explicit construction 
provides that $\psi_Y$ is continuous. 

Suppose also that $Y$ is a Banach analytic space and $\{C_y\}_{y\in Y}$ is 
a holomorphic family of $S$-smooth curves. Assume additionally that
$z_1\equiv z_2$, \ie that $W_{r,R,\rho}$ is a usual (not distorted) cylinder
$A_{r,R} \times \Delta(\rho)$. Let $\psi^-_Y: Y \to \Gamma_{S^-}(A_{r,r'} 
\times \Delta(\rho), \calo)$ and $\psi^+_Y: Y \to \Gamma_{S^+}(A_{R',R}\times
\Delta(\rho), \calo)$ be holomorphic maps, inducing corresponding local 
determinig functions $F^\pm(z,w,y)$ for the family $\{C_y\}$ in $A_{r,r'} 
\times \Delta(\rho)\times Y$ and $A_{R',R}\times \Delta(\rho) \times Y$ 
respectively, see \definition{2.2}. 

\lemma{?2} implies that there exist holomorphic maps $\phi^-_i: Y \to 
\Gamma_{S^-}(A_{r,r'}, \calo)$ such that the Weierstra\3 polynomial $w^n + 
\sum_{i=1}^n f^-_i(z,y) w^{n-i}$ with coefficients\break $f^-_i(z,y)\deff 
\phi^-_i(y)(z)$ is a local 
detemining function for the family $\{C_y\}$ in $A_{r,r'} \times \Delta(\rho) 
\times Y$. Indeed, the map $\phi^-\deff (\phi^-_1, \ldots, \phi^-_n)$ is
obtained as a composition of $\psi^-_Y$ with the map $\Psi$ from \lemma{?2}.
Repeating the same argumentation for $\psi^+_Y$, we obtain a holomorphic map 
$\phi^+ = (\phi^+_1, \ldots, \phi^+_n):Y \to \Gamma_{S^+}(A_{R',R}, \calo^n)$
with similar properties. 

The condition that both $w^n + \sum_{i=1}^n f^\pm_i(z,y) w^{n-i}$ are local
determining functions for the same holomorphic family $\{C_y\}$ means that
the map $(\phi^-, \phi^+): Y \to \Gamma_{S^-}(A_{r,r'}, \calo^n) \times
\Gamma_{S^+}(A_{R',R}, \calo^n)$ takes values in the subset consisting of
the tuples $(f^+,f^-)= ((f^+_1, \ldots,f^+_n),(f^-_1, \ldots,f^-_n))  
\in \Gamma_{S^-}(A_{r,r'}, \calo^n) \times \Gamma_{S^+}
(A_{R',R}, \calo^n)$ which are restrictions onto $A_{r,r'}$ and $A_{R',R}$ of 
some holomorphic function $f=(f_1,\ldots, f_n) \in \Gamma_S(A_{r,R}, \calo^n)$.
This implies that $(\phi^-, \phi^+)$ takes value in $\Gamma_S(A_{r,R}, \calo^n)
\subset \Gamma_{S^-}(A_{r,r'}, \allowbreak
\calo^n) \oplus \Gamma_{S^+}(A_{R',R},
\calo^n)$. Thus any holomorphic family $\{C_y\}_{y\in Y}$ in $A_{r,R} \times 
\Delta(\rho)$ of curves satisfying condition $\supp(C_y) \subset A_{r,R}  
\times  \Delta(\rho')$ with some $\rho' <\rho$ of given degree $n$ is defined 
by a holomorphic map $\phi_Y:Y \to \Gamma_S(A_{r,R}, \calo^n)$.

\smallskip
Now let us return to $W_{r,R,\rho}$ of the general type satisfying the 
conditions of the lemma. Note that all the above constructions of the proof 
respect holomorphic structure, in particular, they can be interpreted as
holomorphic maps between corresponding Banach manifolds. This implies 
that the statement of the lemma about {\sl holomorphic} families 
$\{C_y\}_{y\in Y}$ of $S$-smooth curves is valid. 
\qed

\medskip\newgrph{}
It follows from the above that a small ball in $\Gamma_S(A_{r,R}, \calo^n)$
is a local chart for the space $\calm^{(n)}_S(W_{r,R,\rho})$. 
In this situation an invariant description is also possible.

\newlemma{2.2} The tangent space $T_C\calm_S^{(n)}(W_{r,R,\rho})$ at $C= n\cdot
A_{r,R}$ is canonically isomorphic to $\Gamma_S(A_{r,R}, \caln_C)$.
Formula $\eqq{theta}$ also remains valid.
\endstate

\smallskip\noindent
\bf Proof \rm is identical to that for \lemma{?3}. \qed

\bigskip\bigskip
\newsection{S3}{Globalization}

\nobreak
\newgrph{3.0}
Let $U$ be an open set in a smooth complex surface $X$ and $\calpi^{(1,1)}
(U)$ be the set of all curves in $U$. One can regard $\calpi^{(1,1)} (U)$
as the base of the ``universal'' (weakly continuous) family of curves in $U$. 
However, the weak topology of currents in $\calpi^{(1,1)} (U)$ is not 
convenient to deal with. 
As we have seen in the previous sections, it is more useful to describe
(a family of) curves by appropriate determining functions. Here we shall show
that every continuous family of curves in $U$ can be locally represented
as a continuous deformation of determining functions.

It is enough to consider the situation when $U =\isig \times \Delta$ with 
$\isig$ a smooth complex curve. Let $z$ be a (local) coordinate on $\isig$ and
$w$ a standard one on $\Delta$. Fix a relatively compact subcurve $\isig'
\Subset \isig$ with a smooth boundary and a smoothness $S$ in $\isig'$. Thus 
for every neighbourhood $\iom$ of $\barr\isig{}'$ in $\isig$ the restriction
map $\Gamma (\iom, \calo) \to \Gamma(\isig', \calo)$ takes values in
$\Gamma_S(\isig', \calo) $ and is continuous w.r.t.\ usual Fr\'echet topology
in $\Gamma (\iom, \calo)$.

\newlemma{3.1l} Let $0<r<R<1$ and let $C_0$ be a curve in $U$ whose
restriction $C_0 \cap V_R$, $V_R\deff \isig' \times \Delta(R)$, is $S$-smooth
and lies in $V_r=\isig' \times \Delta(r)$. 
Suppose also that $C_0$ does not contain components of the
form $\{z\} \times \Delta$ with $z\in \barr\isig{}'$. Then there exists a
neighbourhood $\calu^{(n)}$ of $C_0$ in $\calpi^{(1,1)}(U)$ with the
following properties:

\sli For every $C\in \calu^{(n)}$ the restriction $C\cap V_r$ is $S$-smooth
in $V_r$ and has degree $\deg (C \cap V)= \deg(C_0 \cap V)\ddef n$.

\slii The induced map $\kappa: \calu^{(n)} \to \calm_S^{(n)}(V_r)$ is
continuous with respect to the weak topology in $\calpi^{(1,1)}(U)$.

\sliii If $\{C_y\}_{y\in Y}$ is a holomorphic family of curves in $U$ with 
$C_y \in \calu^{(n)}$ for all $y\in Y$, then the induced map $\phi_Y: Y\to 
\calm_S^{(n)}(V_r)$ is holomorphic.
\endstate

\proof{3.1lp} \comment Let $\calu$ be the set of those $C\in \calpi^{(1,1)}
(U)$, such that $|C|\cap \bigl( \barr\isig{}' \times (\barr\Delta(R) \bs
\Delta(r)\bigr) = \emptyset$. \endcomment
Define
$$
\calu \deff \{\, C\in \calpi^{(1,1)} (U) \;:\;
|C|\cap \bigl( \barr\isig{}' \times (\barr\Delta(R) \bs \Delta(r)\bigr) =
\emptyset \,\}.
$$
Since the weak convergence $C_i \longto C$ of currents in $\calpi^{(1,1)}(U)$
implies the Hausdorf convergence of supports $|C_i| \longto |C|$, the set
$\calu$ is open in $\calpi^{(1,1)}(U)$.

Let $\{C_i\}$ be a sequence in $\calu$ which converges to $C\in \calu$ with
$\deg (C_i\cap V_r)=n$. Then \lemma{1.2} implies that $\deg(C\cap V_r)
= n$. Thus $\calu$ is a disjoint union of components $\calu^{(n)} \deff \{\,
C\in \calu \,:\, \deg(C\cap V_r)=n \,\}$ which are open in $\calpi^{(1,1)}
(U)$. Take any $C\in \calu^{(n)}$. Since $U$ is Stein and $\sfh^2(U,\zz)=0$,
there exists $F\in \Gamma(U, \calo)$ such that $C = \div(F)$. But then
$F\ogran_{V_r} \in \Gamma_S(V_r, \calo)$ and this proves the
$S$-smoothness of $C\cap V_r$.

\smallskip
Since $\calpi^{(1,1)}(U)$ is a subset of a space of distributions, its
topology is {\sl sequencial}. This means that a set $A\subset \calpi^{(1,1)}
(U)$ is closed \iff for any sequence $\{C_i\}\subset A$ which
converges to $C=\lim C_i \in \calpi^{(1,1)}(U)$ the limit point $C$ belongs
to $A$. In particular, a map $\kappa: \calu^{(n)} \to \calm_S^{(n)}(V_r)$ is
continuous \iff the image of every convergent sequence is also a convergent
sequence.

So let $\{C_i\}$ be a sequence in $\calu^{(n)}$ converging to $C\in
\calu^{(n)}$. Then there exists a neighbourhood $\iom$ of $\barr\isig'$ in
$\isig$ such that $|C_i| \cap \bigl(\iom \times (\barr\Delta(R) \bs\Delta(r)
)\bigr) = \emptyset$ for every $i>\!\!>1$. It follows that every restricted 
curve $C_i \cap \iom \times \Delta(r)$ is a zero divisor of a uniquely defined
Weierstra\3 polynomial $P_{f_i}$ of degree $n$ with coefficients $f_i=
(f_{i\,1}, \ldots ,f_{i\,n}) \in \Gamma_{L^\infty}(\iom, \calo^n)$. Moreover,
the coefficients $f_i$ are $L^\infty$-bounded uniformly in $i$. Consequently,
$f_i$ weakly converge in $\iom$ to the coefficients $g=(g_1,\ldots,g_n) \in
\Gamma_{ L^\infty}(\iom, \calo^n)$ of the Weierstra\3 polynomial $P_g$ of the
curve $C \cap \iom \times \Delta(r)$.

By the hypotheses of the lemma, the restrictions of $f_i$ onto $\isig'$ are
$S$-smooth and converge to $g\ogran_{\isig'}$ with respect to the norm
topology in $\Gamma_S(\isig', \calo^n)$. This shows that the map $\kappa:
\calu^{(n)} \to \calm^{(n)}(V_r)$ is continuous.

\smallskip
Now let $\{C_y\}_{y\in Y}$ be a holomorphic family of curves in $U$ with all
$C_y \in \calu^{(n)}$. Fix some $y_0\in Y$. Then there exist a neighbourhood
$Y_0$ of $y_0\in Y$ and a neighbourhood $\iom$ of $\barr\isig'$ in $\isig$
such that $|C_y| \cap \bigl(\iom \times (\barr\Delta(R) \bs \Delta(r))\bigr)=
\emptyset$ for every $y\in Y_0$. Take $z^*\in \barr\iom$ and $y^*\in Y_0$, and
consider the set $(\{z^*\} \times \Delta(R)) \cap |C_{y^*}|$. By the
construction, it consists of finitely many points $x_1,\ldots,x_k$. For every
$x_i$ an appropriate multiplicity $m_i$ is defined such that
$\sum_{i=1}^k m_i = n$.

By the definition of a holomorphic family of curves, in some neighbourhood
$W_i \subset U\times Y_0$ of every $(x_i, y^*) $ a holomorphic function
$F_i(z,w;y)\in \Gamma(W_i, \calo_{U\times Y})$ is defined such that
$\div( F_i (\cdot; y)) = C_y \cap W_i$. As in the proof of \lemma{?2}, we can
constract local determinig functions $P_i(z,w;y)\in \Gamma(W_i, \calo_{U
\times Y})$ for $C_y \cap W_i$ which are {\sl polynomial} in $w$,  $P_i
(z,w;y) = w^{m_i} + \sum_{j=1}^{m_i} f_{ij}(z,y)\,w^{m_i-j}$. The product
$P(z,w;y) \deff \prod_{i=1}^k P_i(z,w;y)$ is the Weierstra\3 polynomial of
$C_y$. This shows that in a neighbourhood of $(z^*, y^*)\in \isig \times Y$
the coefficients $f(z;y) = (f_1(z;y) ,\ldots, f_n(z;y) )\in \Gamma_{L^\infty}
(\iom, \calo^n)$ depend holomorphically (\ie analytically) in both variables
$z$ and $y$. It follows that the induced map $\phi_Y: Y\to \Gamma
_{L^\infty} (\iom, \calo^n)$, sending $y\in Y$ into the coefficients of the
Weierstra\3 polynomial of $C_y \cap (\iom \times \Delta(r))$ is holomorphic.
To finish the proof, we apply the restriction map $\Gamma 
_{L^\infty} (\iom, \calo^n) \to \Gamma _S (\isig', \calo^n)$.
\qed

\smallskip
\newgrph{Rem-?} {\bf Remark.} In fact, we have constructed a {\sl morphism}
$\bfphi_Y: Y\to \Gamma_S (\isig', \calo^n)$.

\medskip
\newdefinition{3.1d} Let $U$ be an open set in a smooth complex surface
$X$ and $C$ a curve in $U$. Suppose that there exists a finite collection
$\{U_i\}_{i=1}^N$ of open subsets of $U$ which satisfies the following
properties:

\sli Every $U_i$ is a product $U_i = \isig_i \times\Delta$ with $\isig_i$
being an annulus $A_{r_i,1}$. Moreover, in some neighbourhood $\wt U_i$ of
the closure $\barr U_i$ there exists a holomorphic function $z_i$ whose
restriction on $\isig_i$ coincides with the standard coordinate $z$ on
$A_{r_i,1} = \{\,z \,:\, r_i<|z|<1 \,\}$.

\slii $U \cap \wt U_i = \{ x\in U \;:\;  |z_i(x)| <1 \}$.

\sliii For every $U_i$ there is a fixed smoothness $S_i$ such that $C \cap
U_i$ is a zero divisor of a Weierstra\3 polynomial $P_{g_i}$ of a degree
$n_i$ with $S_i$-smooth coefficients $g_i= (g_{i\,1}, \ldots g_{i\,n_i})
\in \Gamma_{S_i}(\isig_i, \calo^{n_i})$.

\sliv Distinct $U_i$ are disjoint and $|C|\bs \cup_{i=1}^N U_i$ is compact in
$U$.

Then we say that $S\deff \{\, (U_i, z_i, S_i) \,\}$ is a {\sl smoothness in 
$U$} and $C$ is {\sl an $S$-smooth curve in $U$}. A family $\{C_y\}_{y\in Y}$ 
of $S$-smooth curves in $U$ is called {\sl continuous
{\rm (resp.} holomorphic}) \iff $Y$ is a topological (resp.\ complex) space,
$\{C_y\}_{y\in Y}$ is a continuous family of curves in $U$, and every
restricted family $\{C_y \cap U_i\}_{y\in Y}$ is induced by a continuous
(resp.\ holomorphic) map $F_i :Y\to \Gamma_{S_i} (\isig_i, \calo^{n_i})$.

A section $f$ of the structure sheaf $\calo_C$ (resp.\ the normal sheaf
$\caln_C$) is called {\sl $S$-smooth}, $f\in \Gamma_S (C, \calo_C)$ (resp.\
$f\in \Gamma_S (C, \caln_C)$), \iff for every $U_i$ the restriction
$f\ogran_{U_i}$ is $S_i$-smooth. An $S$-smooth curve $C$ is called {\sl
extendible} \iff there exists an (abstract) holomorphic curve $\wt C$ (\ie
a complex analytic space of pure dimension 1) and an open embedding $C
\hookrightarrow \wt C$, such that $|C|$ is relatively compact in $|\wt C|$,
$|C| \Subset |\wt C|$, $\calo_{\wt C}\ogran_C =\calo_C$, and such that the
restriction map $\Gamma(C, \calo_{\wt C}) \longto \Gamma(C, \calo_C)$ takes
values in $\Gamma_S (C, \calo_C)$.

\medskip
\newtheorem{3.1t} Let $X$ be a smooth complex surface, $U\subset X$ an open
subset,
$S$ a smoothness in $U$, and $C$ an $S$-smooth curve in $U$. Suppose that $C$
is extendible. Then there exists a ball $B\subset \Gamma_S (C, \caln_C)$ and
a holomorphic map $\Phi:B \to \sfh^1(C,\caln_C)$ with $\Phi(0) =0$ and
$d\Phi(0)=0$ such that the set $Z\deff \Phi\inv(0)$ is

{\sl a)} a Banach analytic set of finite codimension in $B$ and

{\sl b)} the base of a holomorphic family $\calc= \{\calc_z\}$ of $S$-smooth
curves in $U$ with $\calc_0 =C$ which possesses the following universality
property:

\itemitem{} {\sl For every continuous (resp.\ holomorphic) family $\{C_y\}_{
y\in Y}$ of $S$-smooth curves in $U$ with $C_{y_0} =C$ there exists a
neighbourhood $Y'$ of $y_0$ in $Y$ and a continuous (resp.\ holomorphic) map
$\Psi_Y: Y' \to Z$ with $\Psi_Y(y_0)=0$ and $\calc_{\Psi_Y(y)} = C_y$}.

\rm\medskip
Denote by $\calm_S(U)$ the set of $S$-smooth curves in $U$. Due to
\definition{3.1d} this is a subset of $\prod_i \Gamma_{S_i} (\isig_i,
\calo^{n_i})$ with the induced topology. Thus \theorem{3.1t} provides that in
a neighbourhood of an extendible curve $\calm_S(U)$ has the natural structure
of a Banach analytic space of finite type and that $Z$ is a local chart for
$\calm_S(U)$ at $C$.
We call $\calm_S(U)$ the {\sl moduli space of $S$-smooth curves in $U$}.
\endstate

\proof{3.1tp} Let $\{U_i\}$ be as in \definition{3.1d}. We construct a special
covering $\{V_i\}$ of $|C|$ in $U$ which satisfy the following conditions:

{\sl i$'$)} Every $V_i$ is biholomorphic to $\isig_i \times \Delta$ for some
smooth complex curve $\isig_i$ with a boundary $\d \isig_i$ consisting of
finitely many smooth circles $\gamma_{ij}$, $\d\isig_i =\sqcup_j
\gamma_{ij}$.

{\sl ii$'$)} If $i\le N$, then $V_i = A_{r'_i,1} \times\Delta \subset U_i$ for
some $r_i \le r'_i <1$.

{\sl iii$'$)} With respect to the isomorphism $V_i \cong \isig_i \times
\Delta$, the restricted curve $C\cap V_i$ is a divisor of a Weierstra\3
polynomial $P_i$. Moreover, for every $i>N$ there is fixed a smoothness
$S_i$ on $\isig_i$ and the coefficients of $P_i$ are $S_i$-smooth.

{\sl iv$'$)} An intersection $V_i \cap V_j$ is either empty or is
biholomorphic to a distorted cylinder $W_{ij} \deff W_{r_{ij},R_{ij},
\rho_{ij}}$ with corresponding holomorphic coordinates $z'_{ij}$,
$z''_{ij}$, and $w_{ij}$. In the latter case $|C|\cap V_i \cap V_j$ is a
(non-empty) annulus $\isig_{ij} = A_{r_{ij}, R_{ij}}= \div(w_{ij})$ and
$C_{ij} \deff C\cap V_i \cap V_j =n_{ij} \cdot A_{r_{ij}, R_{ij}}$.

{\sl v$'$)} If $\gamma_{ij}$ is a boundary component of $\isig_i$ with
$i>N$, then $\gamma_{ij} \subset V_j$.

\smallskip
The construction of $\{V_i\}$ can be realised as follows: First, for
every $i\le N$ we find $r'_i$ with $r_i \le r'_i <1$ such that $|C|$ is a
smooth analytic set in a neighbourhood of $|C| \cap (\{ x\in \isig_i \,:\,
|z(x)|=r'_i \} \times \Delta)$. Next, we consider the singular points of
$|C|\bs (\cup_{i=1}^N V_i)$ and find an appropriate neighbourhood $V_i$,
$i=N+1, \ldots, N_1$, of every such a point, so that $V_i$ and $V_j$ are
disjoint for $1\le i,j\le N_1$. Then the set $|C|\bs (\cup_{i=1}^{N_1} V_i)$
can be covered by finitely many smooth complex non-closed curves $C'_k$ with
a smooth boundary which we enumerate by $k=N_1 +1,\ldots N_2$.

For any $C'_k$ we fix a neighbourhood $V'_k$ of a closure $\barr C{}'_k$
such that $|C| \cap V'_k$ is also smooth with a smooth boundary and
$\sfh^2(V'_k, \zz)=0$. In particular, the (holomorphic) line bundle $L_{C \cap
V'_k}$, corresponding to a divisor $|C| \cap V'_k$, is topologically trivial.
Due to a result of Siu [Siu], the set $|C| \cap V'_k$ admits a Stein
neighbourhood $V''_k \subset V'_k$. The condition of topological triviality
of $L_{C \cap V'_k}$ provides the existence of a holomorphic function $w_k
\in \Gamma(V''_k, \calo)$ such that $|C| \cap V''_k= \div(w_k)$.

We may assume that $|C| \cap V''_k$ is biholomorphic to a subdomain of the
complex plane $\cc$. Let $z_k$ be a holomorphic function on $|C| \cap V''_k$
which corresponds to a standard coordinate on $\cc$. Since $V''_k$ is Stein,
we can extend $z_k$ to a holomorphic function in $V''_k$. Now one can see
that,
choosing appropriate $\isig_k \subset |C| \cap V''_k$, $k=N_1 +1,\ldots
N_2$, and setting $V_k\deff \{\, x\in V''_k \,:\, z_k(x) \in \isig_k |w_k(x)|
<r_k \,\}$,
it is possible to obtain the desired covering $\{V_i\}$ with $i=1,\ldots,N_2$.

\smallskip
\newgrph{G3}
Due to the construction of $V_i = \isig_i \times \Delta$, the boundary
components $\{\gamma_{ij}\}$ of $\isig_i$ are naturally separated into two
groups which consist respectively of ``inner'' components lying in $U$ and
``outer'' components lying on $\d U$. It is easy to see that the property of
a curve $C$ to be $S$-smooth in $U$ is independent of the choice of inner
smoothness classes $S_{ij}$ which correspond to inner components
$\gamma_{ij}$.
Thus without loss of generality we may assume that all inner smoothnesses
classes $S_{ij}$ are Hilbert, \ie the corresponding spaces $\Gamma_{S_i}
(\Delta, \calo)$ are Hilbert spaces. For example, one can take all $S_i$ to
be some Sobolev smoothness class $L^{k,2}$.

\smallskip
For any index pair $(i,j)$ with nonempty $W_{ij} = V_i \cap V_j$ we denote by
$n_{ij}$ the multiplicity of $C_{ij} =C\cap W_{ij}$. Note that for such
$(i,j)$ the smoothnesses $S_i$ and $S_j$ in $V_i = \isig_i \times \Delta$ and
$V_j= \isig_j \times \Delta$ induce the smoothness $S_{ij}$ on $\isig_{ij}$
with the continuous projections $\Gamma_{S_i}(V_i,\calo) \to
\Gamma_{S_{ij}}(W_{ij}, \calo)$ and $\Gamma_{S_j}(V_j, \calo) \to
\Gamma_{S_{ij}} (W_{ij}, \calo)$. We fix sufficiently small balls $B_{ij}
\subset \Gamma_{S_{ij}} (\isig_{ij}, \calo^{n_{ij}})$ which parameterise
$S_{ij}$-smooth curves in $W_{ij}$ which are sufficiently close to $C\cap
W_{ij}$.

\smallskip
Now fix some $V_i$. Then the restricted curve $C \cap V_i$ is a zero divisor
of a uniquely defined Weierstra\3 polynomial $P_i$ of the degree $n_i$ with
$S_i$-smooth coefficients $g_i=(g_{i\,1},\ldots) \in \Gamma_{S_i} (\isig_i,
\calo^{n_i})$. Fix a sufficiently small ball $B_i \deff B(g_i, a_i) \subset
\Gamma_{S_i} (\isig_i, \calo^{n_i})$ centered at $g_i$. If the radius $a_i$
of $B_i$ is chosen sufficiently small, then, for every $j$ such that $W_{ij}
\not= \emptyset$ and for every $f\in B_i$, the restricted curve $\div(P_f)
\cap W_{ij}$ is a zero divisor of a uniquely defined distorted Weierstra\3
polynomial $\wt P_g$ of the degree $n_{ij}$ with $S_{ij}$-smooth coefficients
$g\in \Gamma_{S_{ij}} (\isig_{ij}, \calo^{n_{ij}})$.

This defines a map $\phi_{ij}: B_i \to B_{ij}$ which is holomorphic (see
\lemma{2.1}). We may assume that the image of $\phi_{ij}$ lies in the ball
$\half B_{ij}$. 
Consider the holomorphic map
$$
\wt\Phi: \prod_i B_i \lrar \prod_{i<j} B_{ij},
\qquad
(\wt\Phi(f))_{ij}\deff \phi_{ij}(f_i) -
\phi_{ji}(f_j)
$$
The product $\prod_i B_i$ parametrised the space of tuples $(C_i)$, where every
$C_i\deff \div(f_i)$ is an $S_i$-smooth curve in $V_i$ which is
 sufficiently close to
$C\cap V_i$. Two such ``local deformations'' $C_i$ and $C_j$ coincide
exactly when $\phi_{ij}(f_i) = \phi_{ji}(f_j)$. It follows that the
analytic set $\wt Z\deff \wt\Phi{}\inv(0)\subset \prod_i B_i$ satisfies
property \slii of the theorem.

\medskip
Due to {\sl Lemmas \lemmano{?3}, \lemmano{2.2}} and \definition{GsN}, the
tangent space to $\prod_i B_i$ at $g=(g_i)$ is isomorphic to 
$ 
\sum_i\Gamma_{S_i}(\isig_i, \calo^{n_i}) = 
\sum_i \Gamma_{S_i}(C\cap V_i, \caln_C),
$
whereas the tangent space to $\prod_{i<j} B_{ij}$ at $0$ is isomorphic to
$ \sum_{i<j}\Gamma_{S_{ij}} (\isig_{ij}, \calo^{n_{ij}}) = \sum_{i<j}
\Gamma_{S_{ij}}(C\cap V_{ij}, \caln_C) $. Formula $\eqq{theta}$ implies
that the differential $d\wt\Phi_g$ of $\wt\Phi$ at $g$ coincides with the
\v{C}ech coboundary operator
$$
d\wt\Phi_g= \delta_S: \sum_i \Gamma_{S_i}(C\cap V_i, \caln_C) \lrar
\sum_{i<j} \Gamma_{S_{ij}}(C\cap V_{ij}, \caln_C),
\quad (\delta(f_i))_{ij}= f_i\ogran_{W_{ij}} - f_j\ogran_{W_{ij}}.
$$

The key point of the proof is that for an extendible curve $C$ the
operator $\delta$ has a closed image and splits. To show this we fix some
$i\le N$ so that $V_i= \isig_i \times\Delta$ touches the boundary $\d U$.
Let $\gamma_{i1}$ be a boundary component of $\isig_i$ lying on $\d U$ and
$S_{i1}$ the corresponding smoothness class. Then there exists curves $C'_i
\subset C''_i \subset \wt C$ such that:

{\sl a)} $|C'_i| \cap |C|= |C''_i| \cap |C|= |C_i|$ so that both $C'_i$ and
$C''_i$ are extensions of $C_i$ ``outwards" from $|C|$.

{\sl b)} $|C'_i|$ is relatively compact in $|\wt C|$ and the ``outer'' part
of the (topological) boundary of $|C'_i|$ lying outside $|C|$ is smooth and
consists of finitely many circles $\gamma'_{ij}$ which lie in
$|C''_i|$.

{\sl c)} If $f\in \Gamma(C'_i, \calo_{\wt C})$, then $f\ogran_{C'_i,}$ is
$S_{i1}$-smooth at the ``outer'' part of the boundary of $|C_i|$ which lie on
$\d U$.

We repeat this construction for every $1\le i\le N$ and set $C'_i = C''_i =
C_i$ and so on for $i>N$. Set $C'\deff C \cup_i C'_i$ and $C''\deff C \cup_i
C''_i$. These are complex curves. The boundary $\d C'$ of $C'$ consists of
smooth circles $\gamma'_{ij}$. Since the restriction $\caln_C\ogran_{C_i}$ is
trivial, we can extend $\caln_C$ to a rank 1 locally free $\calo_{C''}$-module
$\caln_{C''}$ with $\caln_{C''}\ogran_{C''_i}$ trivial. For any component
$\gamma'_{ij}$ of the ``outer'' part of the boundary of $|C'_i|$ we fix a
Hilbert smoothness class $S'_{ij}$. This defines the Hilbert space
$\Gamma_{S'_i}(C'_i, \caln_{C'})$ and a (continuous) restriction map
$\Gamma_{S'_i}(C'_i, \caln_{C'}) \to \Gamma_{S_i}(C_i, \caln_C)$.

Consider the induced \v{C}ech coboundary operators
$$
\delta': \sum_i \Gamma_{S'_i}(C'_i, \caln_{C'}) \lrar
\sum_{i<j} \Gamma_{S_{ij}}(C_{ij}, \caln_C),
\qquad
(\delta'(f_i))_{ij}\deff f_i\ogran_{C_{ij}} - f_j\ogran_{C_{ij}},
\neqq{delta'}
$$
and
$$
\delta'': \sum_i \Gamma(C''_i, \caln_{C''}) \lrar
\sum_{i<j} \Gamma(C_{ij}, \caln_C),
\qquad
(\delta''(f_i))_{ij}\deff f_i\ogran_{C_{ij}} - f_j\ogran_{C_{ij}}.
\neqq{delta''}
$$

By the construction, all $C''_i$ are Stein spaces. Thus \eqq{delta''} is an
acyclic \v{C}ech resolvent for $\caln_{C''}$. Consequently, $\ker(\delta'') =
\sfh^0(C'', \caln_{C''})= \Gamma(C'', \caln_{C''})$ and $\coker(\delta'') =
\sfh^1(C'', \caln_{C''})$. We note the canonical isomorphisms
$$
\sfh^1(C'', \caln_{C''}) \cong \sfh^1(C', \caln_{C'})
\cong\sfh^1(C, \caln_C).
$$
These are finite dimensional spaces. Denote by $p$ the composition
$$
\sum_{i<j} \Gamma_{S_{ij}}(C_{ij}, \caln_C)
\hookrightarrow \sum_{i<j} \Gamma(C_{ij}, \caln_C) \twoheadrightarrow
\sfh^1(C'', \caln_{C''})
$$
and set $T \deff \ker(p)$.

First we note that $p$ is a surjection onto $\sfh^1(C'', \caln_{C''})$.
For this observe that one can find an acyclic \v{C}ech resolvent
$$
\hat\delta: \sum_i \Gamma(\wh C_i, \caln_{C''}) \lrar
\sum_{i<j} \Gamma(\wh C_{ij}, \caln_C),
\qquad
(\hat\delta(f_i))_{ij} \deff
f_i\ogran_{\wh C_{ij}} - f_j\ogran_{\wh C_{ij}},
\neqq{hdel}
$$
for $\caln_{C''}$ with $C_{ij} \Subset \wh C_{ij}$. Then every $[h]\in
\sfh^1(C'', \caln_{C''})$ can be represented by $h=(h_{ij})$ with
$h_{ij} \in \Gamma(\wh C_{ij}, \caln_C)$ so that the restriction
gives $h\in \sum_{i<j} \Gamma_{S_{ij}}(C_{ij}, \caln_C)$.

Now take $h=(h_{ij})\in T$. Since $p(h)=0$, there exists $f''=(f''_i) \in\sum_i
\Gamma(C''_i, \caln_{C''})$ such that $h= \delta''(f)$. Let $f'_i\in
\Gamma(C'_i, \caln_{C'})$ denote the restriction of $f''_i$ onto $C'_i$.

Now in fact $f'_i\in \Gamma_{S'_i}(C'_i, \caln_{C'})$. The
corresponding smoothness of $f'_i$ at the outer component $\gamma'_{ij}$ of
the boundary $\d C'_i$ follows from the fact that $f'_i$ is holomorphic in
a neighbourhood of $\gamma'_{ij}$. Similarly, if $\gamma_{ij}$ is an inner
component of the boundary $\d C'_i$, then $\gamma_{ij}$ lies in some $V_j$.
In this case $f'_j$ is holomorphic in a neighbourhood of $\gamma_{ij}$ and
$f'_i = h_{ij} + f'_j$. Since $h_{ij}$ is $S_{ij}$-smooth at $\gamma_{ij}$,
the same holds for $f'_i$.

This implies that the image of $\delta'$ is $T$ and is of finite
codimension. Consequently $T$ is a closed subspace of $\sum_{i<j}
\Gamma_{S_{ij}}(C_{ij}, \caln_C)$. Since all the smoothnesses $S_i$ are
Hilbert, $\sum_i \Gamma_{S'_i}(C'_i, \caln_{C'})$ is a Hilbert space
and $\ker(\delta')$ admits a complement. Therefore there exists a
splitting operator $\sigma': \Gamma_{S_{ij}}(C_{ij}, \caln_C) \to \sum_i
\Gamma_{S'_i}(C'_i, \caln_{C'})$ such that for every $h\in T$ holds
$\delta'(\sigma'(h))=h$.

Let $\sigma: \Gamma_{S_{ij}}(C_{ij}, \caln_C) \to \sum_i \Gamma_{S'_i}(C'_i,
\caln_{C'})$ denote the composition of $\sigma'$ with the restriction map
$\sum_i \Gamma_{S'_i}(C'_i, \caln_{C'}) \to \sum_i \Gamma_{S_i}(C_i,
\caln_C)$. Then again $\delta(\sigma(h))=h$.

\medskip
Recall that $g=(g_i)\in \prod_i B_i$ is a tuple of coefficients Weierstra\3 
polynomials parametrising our curve $C$ and $\delta$ is the differential 
$d\wt\Phi_g$. Denote by $\wt\Phi_T$ the composition of $\wt\Phi$ with the 
orthogonal projection on $T$ and by $\Phi$ the composition of $\wt\Phi$ 
with the projection onto $\coker d\wt\Phi_g = \sfh^1(C, \caln_C)$.

The implicit function theorem implies that the set $Z_1\deff \wt\Phi_T
\inv(0)$ is a complex Banach submanifold of $\prod_i B_i$ containing
$g$ with the tangent space $T_gZ_1$ canonically isomorphic to $\ker(d\wt\Phi_g)
= \sfh^0_S(C, \caln_C)$. One can easily see that the set $Z_2= \wt \Phi
\inv (0)$ is an analytic subset of $Z_1$ defined by the equation $Z_2 =
\{\,y\in Z_1:  \Phi(y) =0 \,\}$. Take a neighbourhood $B$ of $g$ in $Z_1$
biholomorphic to a small ball in $\sfh^0_S(C, \caln_C)$ and set $Z\deff Z_2
\cap B$.  Then $B$, $\Phi$, and $Z$ satisfy the condition of the theorem.
\qed

\bigskip\bigskip
\newsection{S4}{Applications and related questions.}

\nobreak
\newgrph{4.0} 
\bf Proof of the Main Theorem. \rm Let $X$ be a smooth complex surface, $C^*$ 
a curve in $X$, and $K\Subset |C^*|$ some compact subset. Repeating the
constructions of the proof of \theorem{3.1t}, one can find an open
neighbourhood $U \subset X$ of $K$ and appropriate collections $\{U_i\}$
and $\{V_i\}$ of open sets in $U$ such that $C^* \cap U$ is $S$-smooth with
respect to an appropriate smoothness $S$ in $U$. Furthermore, $\{U_i\}$ and
$\{V_i\}$ can be chosen to
satisfy the conditions \sli--\sliv of the \definition{3.1d} and
{\sl i$'$)--v$'$)} of the proof of \theorem{3.1t}. Then the statement of
the {\sl Main Theorem} follows immediately \theorem{3.1t}.
\qed

\medskip
\newgrph{Ram}
It is important to underline the fact that Banach
analytic sets of finite type, in contrast to general Banach analytic sets,
have a simple structure (see [Ra], Chapter II, \S\.3).
Namely, let $B$ be a ball in a Banach space $E$ and let $Z$ be an analytic
subset of $B$, which is defined by finitely many holomorphic
function and contains $0\in E$. Then in a neighbourhood $V\subset B$ of $0$ 
the set $Z$ has finitely many components, $Z\cap V = \cup_{i=1}^N Z_i$, 
each of them is irreducible at $0$ and is
defined by finitely many holomorphic functions too. Further, for every such
component $Z_i$ of $Z$ at $0$ there exist a closed (Banach) subspace $E_i
\subset E$ of finite codimension and a (linear) projection $\pi_i: E\to E_i$
such that in some smaller neighbourhood $V_i\subset V$ of $0$ the restricted
projection $\pi_i: Z_i\cap V_i \to \pi_i(V_i)$ is a proper branched analytic
covering with $\pi_i\inv(0) \cap Z_i =\{0\}$.

\medskip
As a corollary of the {\sl Main Theorem} we obtain the following statements.

\newproposition{Cn} Let $X$, $C^*$, and $U$ be as in the {\sl Main Theorem}.
Suppose that $\{C_n\}$ is a sequence of curves in $X$ converging (weakly) to
$C^*$. Then for any $n$ which is sufficiently big there exists
a holomorphic family
$\{C_\lambda\} _{\lambda \in \Delta}$ of curves in $U$ which is parametrised
by a disk $\Delta \subset \cc$ and contains both $C_n\cap U$ and $C^*\cap U$.
\endstate

\proof{Cnp} Denote $C\deff C^*\cap U$. Let $B \subset \Gamma_S(C, \caln_C)$
be a small ball and $Z \deff \Phi\inv(0)\subset B$ a local chart for $\calm_S
(U)$. Let $Z = \cup Z_i$ be the decomposition of $Z$ into components such
that every $Z_i$ is a proper branched analytic covering of a ball $B_i$ in an
appropriate Banach subspace $E_i\subset \Gamma_S(C, \caln_C)$ with respect to
a projection $\pi_i: \Gamma_S(C, \caln_C) \to E_i$. Then for $n>\!\!>1$ a
curve $C_n \cap U$ is parametrised by a uniquely defined $a_n \in Z$, in
particular, $a_n$ lies in some $Z_i$.

Set $a'_n \deff \pi_i (a_n)$ and let $L_n$ be a complex line in $E_i$ through
$a'_n$ and $0$. Then $\pi_i\inv (L_n) \cap Z_i$ is a complex curve which
consists of finitely many irreducible components, each of which contains
$0$. Consequently, there exists a holomorphic map $f_n: \Delta \to \Gamma_S
(C, \caln_C)$ such that $f_n(\lambda) \in Z_i$, $f_n(0) =0$, and $f_n(
\lambda_n) = a_n$ for some $\lambda_n \in \Delta$. The map $f_n$ defines the
desired one parameter family $\{C_\lambda\}_{\lambda \in \Delta}$ connecting
$C^* \cap U$ and $C_n \cap U$.
\qed

\smallskip
\newgrph{Le} As an application, the {\sl Main Theorem} and \proposition{Cn}
yields a generalization of Levi's continuity principle. For a precise
statement we need the notion of a {\sl meromorphic hull of a domain $W$ in a
complex manifold $X$}. Recall, that this is a maximal Riemann domain $(\wh W,
\pi)$ over $X$ containing $W$, (\ie $\pi: \wh W \to X$ is a locally
biholomorphic map and there is given an inclusion $i: W \to \wh W$ with $\pi
\scirc i= \id_W$), such that every meromorphic function $f$ on $W$ extends to
a function $\hat f$ on $\wh W$. We refer to [Iv] and [IvSh] for details. 
\comment
Usually, one denotes a meromorphic hull simply by $\wh W$, having in mind the
corresponding projection map $\pi$ and an inclusion $i_W$.
\endcomment

 Let $X$ be a complex surface,
$W\subset X$ a domain, $\wh W$ its meromorphic hull, and $\{ \wh C_n\}$ a
sequence of curves in $\wh W$ without multiple components. Suppose that there
exists a domain $U\Subset X$, such that the projected curves $C_n\deff \pi
(\wh C_n)$ are close in $U$ and (weakly) converge to a curve $C_\infty$.
Suppose also the boundary of $C_\infty$ is not empty and lies in $W$. Then the
sequence $\{ \wh C_n\}$ converge to a curve $\wh C_\infty$ with $\pi(\wh
C_\infty)= C_\infty$.
\endstate

\newgrph{CPLr}
\bf Remark. \rm\theorem{CPL} has the following meaning. If $f$ is a
meromorphic function in $W$ which extends meromorphically to a neighbourhood
of every $C_n$ (possibly as a several sheeted function), then it can be
extended in a neighbourhood of $|C_\infty|$. Thus we obtain the
generalization of the classical result of E.\.E.\.Levi which deals with 
the case where $C_\infty$ and every $C_n$ are disks.

\proof{CPLp} Set $K_0\deff |C_\infty| \bs W$. Let $K$ be the union of $K_0$
with those connected components of $|C_\infty| \cap W$ which are relatively
compact in $|C_\infty|$. Then $K$ is compact. According to the {\sl Main
Theorem}, we can find an open set $U_1\subset U$ and a smoothnes $S$ in
$U_1$ such that $K\subset U_1$ and $C_\infty \cap U_1$ is $S$-smooth.

Due to \proposition{Cn}, for any $n>\!\!>1$ there exists a holomorphic family
$\{Z_\lambda\}_{\lambda \in \Delta}$ of curves in $U_1$ which contains
containing $C_n$ and $C_\infty$. Let $\lambda_0$ and $\lambda_\infty$ be the
corresponding parameter values. If $n$ was chosen big enough, the boundary of
every $Z_\lambda$ lies in $W$. In an obvious way, this family defines a
complex space $Z \subset U_1 \times \Delta$ such that every $Z_\lambda$ is
identified with $Z\cap U_1 \times \{\lambda\}$.

Let $f$ be a meromorphic function in $W$ and $\hat f$ its extension on
$\wh W$. Using the technique of [Iv] and [IvSh] for meromorphic extension
along a holomorphic family of curves, one can show that the restriction of
$f$ onto $Z_{\lambda_0}$ extends to a holomorphic function $F$ on the entire
space $Z$ such that $F$ coincides with $f$ in a neighbourhood of boundary of
every $Z_\lambda$. Since $f$ can be any meromorphic function in $W$, this
means that the family $\{Z_\lambda\}_{\lambda \in \Delta}$ can be lifted to
a family $\{\wh Z_\lambda\}_{\lambda \in \Delta}$ of curves in $\wh W$ such
that $\pi(\wh Z_\lambda)= \wh Z_\lambda$ and $\wh Z_\lambda \cap W= Z_\lambda
\cap W$. For further details see [Iv] and [IvSh].

Now it is not difficult to show that the desired curve $\wh C_\infty\subset
\wh W$ can be constructed as $ \wh Z_{\lambda_\infty} \cup (C_\infty \cap W)$.
\qed

\smallskip
\newgrph{hi-dim} Trying to generalise the results of the paper one must 
overcome the following difficulties:

{\sl1)} Considering deformation problem a non-compact complex subspace $Z$ 
in a complex manifold $X$, such that $\dimc Z>1$, one confronts with the 
fact that non-compact components of $Z$ can be non-Stein. Thus an appropriate
cohomology group $H^1(Z, \caln_Z)$ can be infinite-dimensional and even worse
not separated topological vector space.

{\sl2)} We illustrate problem which appears by deformation of non-compact 
cycles by the following example. Let $n\ge2$ and $k\ge2$ be integers. Consider
disc $\Delta_0\deff \Delta \times \{0\} \subset \Delta \times \cc^n$ and 
cycle $Z \deff k\cdot \Delta_0$ in $X \deff \Delta \times \cc^n$. The problem 
of deformation of $Z$ leads to consideration of the space $\hom(\Delta, 
\sym^k \cc^n)$ of holomorphic maps between $\Delta$ and $k$-th symmetric power
of $\cc^n$. The space $\sym^k\cc^n$ is naturally realised as an analytic subset
of some $\cc^N$. This gives a natural inclusion $\hom(\Delta, \sym^k \cc^n)
\subset \hom(\Delta, \cc^N)$. Fixing some smoothness class $S$, \eg $S=
L^\infty$, we obtain the set $\hom_S(\Delta, \sym^k \cc^n)$ of $S$-smooth maps
as an {\sl analytic} subset of the Banach space $\hom_S(\Delta, \cc^N)$. Thus
$\hom_S(\Delta, \sym^k \cc^n)$ is equipped with the natural structure of Banach
analytic space. However, $\hom_S( \Delta, \sym^k \cc^n)$ is has infinite
codimension in $\hom_S(\Delta, \cc^N)$. Moreover, there are infinitely many 
irreducible components of $\hom_S(\Delta, \sym^k \cc^n)$ in $f \equiv0$ (zero 
map).

{\sl3)} One can obtain a statement similar to \theorem{3.1t} considering 
deformation of stable maps from non-compact nodal curves to a given
smooth complex
manifold $X$ of arbitrary dimension. Recall, that an abstract {\sl nodal curve 
$C$} is a complex space of dimension one whose singularities are only ordinary 
double points. Such curves are also called semi-stable. A holomorphic map
$f:C \to X$ is called stabe, if there are only finitely many biholomorphisms
$g: C \to C$ with the property $f\scirc g =f$. In [IvSh] the deformation 
problem of the stable maps $(C,f)$ is considered under the following additional
assumtions. 

\itemitem{$(*)$}
Curves $C$ have finitely many irreducible components, each of 
them being of finite genus and bordered by finitely many smooth circles, and
maps $f:C \to X$ are $L^{1,p}$-smooth up to boundary $\d C$. 

The set of such pairs $(C,f)$ is equipped with {\sl Gromov topology}. The 
following result is proved.

\newgrph{st-map} \sl For given $(C_0, f_0)$ there exist Banach analytic spaces
$\calm$, $\calc$ of finite type and holomorphic maps $\pi:\calc \to \calm$,
$F: \calc \to X$, such that:

\sli for any $\lambda \in \calm$ the fiber $C_\lambda\deff \pi\inv(\lambda)$ 
is a nodal curve and $f_\lambda\deff F\ogran_{C_\lambda}: C_\lambda \to X$ 
is a stable map with the property $(*)$; 

\slii for any stable map $(C,f)$ with the property $(*)$ sufficiently close 
to $(C_0, f_0)$ w.r.t.\ Gromov topology there exist a $\lambda \in \calm$ and 
biholomorphism $\phi: C \to C_\lambda$ such that $f= f_\lambda \scirc \phi$.

\smallskip

\spaceskip=4pt plus3.5pt minus 1.5pt
\xspaceskip=5pt plus4pt minus 2pt
\font\sc=cmcsc10
\newdimen\length
\newdimen\lleftskip
\lleftskip=3.3\parindent
\length=\hsize \advance\length-\lleftskip
\def\entry#1#2#3#4\par{\parshape=2  0pt  \hsize%
\lleftskip \length%
\noindent\hbox to \lleftskip%
{\bf[#1]\hfill}{\sc{#2 }}{\sl{#3}}#4%
\medskip
}

\bigskip\bigskip
\centerline{\bigbf References.}

\medskip

\entry{Ba}{D.\.Barlet.}{Espace analytique r\'eduit des cycles analytiques
complexes d'un espace analytique complexe de dimension finie.}
Lecture notes in mathematics {\bf482}, Springer, 1975.

\entry{Dou}{A.\.Douady.}{Le probl\`em des modules pour les sous-espaces
analytique compacts d'une espace analytique donn\'e.} Ann.\ Inst.\
Fourier, {\bf16}, 1--95,  (1966).

\entry{Ga}{Th.\.Gamelin.}{Uniform algebras.} Prentice-Hall, Englewood Cliffs,
N.J., 1969.

\entry{GrHa}{P.\.Griffiths, J.\.Harris.}{Principles
of algebraic geometry.} Wiley\&Sons, New York, 1978.

\entry{Ha}{R.\.Harvey.}{Holomorphic chains and their boundaries.} Proc.\
Symp.\ Pure Math.\ at Williamstown, 1975, p.309--382.

\entry{Iv}{S.\.Ivashkovich.}{The Hartogs-type extension theorem for
meromorphic maps into compact K\"ahler manifolds.} Invent. Math. {\bf109},
47--54, 1992.


\entry{IvSh}{S.\.Ivashkovich, V.\.Shevchishin.}{Defotmation of non-compact
curves and envelopes of meromorphy of spheres.} To appear in Mathematical 
Sbornik.

\entry{Ki}{J.\.King.}{The currents defined by anylytic varieties.} Acta Math.,
{\bf127}, 185--220, (1971).

\entry{Ra}{J.-P.\.Ramis.}{Sous-ensembles analytiques d'une vari\'et\'e
banachique complexe.} Springer, Berlin, 1970.

\entry{Siu}{Y.-T.\.Siu.}{Every Stein subvariety admits a Stein
neighbourhood.} Invent.\ math., {\bf38}, 89--100, 1976.

\enddocument
\bye